%


\magnification=\magstep1
\def\forces{\parallel\!\!\! -}
\def\Bool{[\![}
\def\Boor{]\!]}


\def\hexnumber#1{\ifcase#1 0\or1\or2\or3\or4\or5\or6\or7\or8\or9\or
	A\or B\or C\or D\or E\or F\fi }

\font\teneuf=eufm10
\font\seveneuf=eufm7
\font\fiveeuf=eufm5
\newfam\euffam
\textfont\euffam=\teneuf
\scriptfont\euffam=\seveneuf
\scriptscriptfont\euffam=\fiveeuf


\font\tenmsx=msam10
\font\sevenmsx=msam7
\font\fivemsx=msam5
\font\tenmsy=msbm10
\font\sevenmsy=msbm7
\font\fivemsy=msbm5
\newfam\msxfam
\newfam\msyfam
\textfont\msxfam=\tenmsx  \scriptfont\msxfam=\sevenmsx
  \scriptscriptfont\msxfam=\fivemsx
\textfont\msyfam=\tenmsy  \scriptfont\msyfam=\sevenmsy
  \scriptscriptfont\msyfam=\fivemsy
\edef\msx{\hexnumber\msxfam}

\mathchardef\upharpoonright="0\msx16
\let\restriction=\upharpoonright
\def\Bbb#1{\tenmsy\fam\msyfam#1}

\def\restrict{{\restriction}}
\def\Veskip{\vskip0.8truecm}
\def\Smallskip{\vskip1.4truecm}
\def\Bigskip{\vskip2.2truecm}
\def\Hoskip{\hskip1.6truecm}

\def\qed{{\vcenter{\hrule height.4pt \hbox{\vrule width.4pt height5pt
 \kern5pt \vrule width.4pt} \hrule height.4pt}}}
\def\notin{{\in}\kern-5.5pt / \kern1pt}

\def\BB{{\Bbb B}}

\def\EE{{\Bbb E}}
\def\FF{{\Bbb F}}
\def\II{{\Bbb I}}

\def\LL{{\Bbb L}}
\def\LA{{\Bbb{LA}}}
\def\MM{{\Bbb M}}

\def\MI{{\Bbb{MI}}}

\def\PP{{\Bbb P}}
\def\QQ{{\Bbb Q}}

\def\sm{{\smallskip}}
\def\ce#1{{\centerline{#1}}}
\def\no{{\noindent}}
\def\la{{\langle}}
\def\ra{{\rangle}}
\def\sub{\subseteq}
\def\tr{{\triangleright}}
\font\small=cmr8 scaled\magstep0

\font\capit=cmcsc10 scaled\magstep0
\font\capitg=cmcsc10 scaled\magstep1

\font\dunhgg=cmr10 scaled\magstep2

\font\bold=cmssdc10 scaled\magstep1
\font\bolds=cmssdc10 scaled\magstep0
\font\slan=cmssdc10 scaled\magstep0
\overfullrule=0pt
\openup1.5\jot

\centerline{}
\Bigskip
\centerline{\dunhgg Combinatorial properties of classical forcing
notions}
\Bigskip
\centerline{\capitg J\"org Brendle\footnote{$^*$}{{\small
The author wishes to thank the MINERVA-foundation
for supporting him}}}
\Smallskip
{\baselineskip=0pt {\small
\noindent 
Department of Mathematics,
Bar--Ilan University,
52900 Ramat--Gan, Israel
\par
and \par
\noindent Mathematisches Institut der Universit\"at T\"ubingen,
Auf der Morgenstelle 10, 72076 T\"ubingen, Germany}} 
\Bigskip
\centerline{\capitg Abstract}
\bigskip
\noindent We investigate the effect of adding a single real
(for various forcing notions adding reals) on cardinal invariants
associated with the continuum. We show: (1) adding an
eventually different or a localization real adjoins a Luzin
set of size continuum and a mad family of size $\omega_1$;
(2) Laver and Mathias forcing collapse the dominating number
to $\omega_1$, and thus two Laver or Mathias reals added
iteratively always force $CH$; (3) Miller's rational perfect
set forcing preserves the axiom $MA(\sigma$--centered).
\vfill\eject
{\bold Introduction}
\Smallskip
In [Ro], Roitman proved 
that adding a Cohen real preserves
the axiom $MA(\sigma$--centered). However,
larger fragments of $MA$ can fail in the
extension. In fact, Woodin [Wo], and
(independently, but later) Cicho\'n and Pawlikowski [CP]
showed that $MA$ for random forcing fails after adjoining a Cohen
real. The latter authors investigated as well the effect of adding
a Cohen or a random real on cardinal invariants in Cicho\'n's
diagram (see [CP] and [Pa 2] --- these are cardinals associated
with the ideals of meager and null sets on the real line as well 
as the unbounding and dominating numbers --- see section 1 
for definitions).
Recently, Judah, Shelah and the present author did the same for
Hechler forcing [BJS]. We continue this line of research for 
other forcing notions adding a single real 
--- in particular for forcings involved in
consistency proofs in Cicho\'n's diagram ---; and the present
work should be looked at as a companion piece to [BJS].
\par
After reviewing the forcing notions as well as the cardinal
invariants we shall be considering in our work in section 1,
we will show in 2.1. that many $ccc$ forcing notions
(including Miller's eventually different reals forcing
[Mi 1, section 5] as well as Bartoszy\'nski's localization
forcing) add a Luzin set of size $2^\omega$. Remember that an
uncountable set of reals $X\sub2^\omega$ is {\it Luzin}
iff it has at most countable intersection with every
meager set. --- We shall then sketch how some changes in
the argument prove that the forcings considered in 2.1.
produce a maximal almost disjoint family of subsets of 
$\omega$ of size $\omega_1$ (Theorem 2.3.). Recall that
$A, B \sub \omega$ are said to be {\it almost disjoint}
({\it a. d.} for short) iff $\vert A \cap B \vert
< \omega$; ${\cal A} \sub [\omega]^\omega$ is an {\it
a. d. family} iff the members of ${\cal A}$ are pairwise
a. d.; and ${\cal A}$ is a {\it mad family} ({\it maximal
almost disjoint family}) iff it is a. d. and maximal with
this property. --- We close section 2 with consequences of
adding an eventually different (localization) real on the
cardinals in Cicho\'n's diagram.
\par
The third section is devoted to Laver and Mathias forcings
(and their versions with ultrafilters). Recall that given reals
$f,g \in\omega^\omega$, $f$ {\it eventually dominates} $g$ iff
$\forall^\infty n \; (f(n) > g(n))$ (where $\forall^\infty
n$ denotes {\it for almost all $n$}; similarly $\exists^\infty
n$ stands for {\it there exist infinitely many $n$}).
${\cal F} \sub \omega^\omega$ is a {\it dominating family}
iff for all $g\in\omega^\omega$ there is $f \in {\cal F}$
eventually dominating $g$. We shall prove in 3.1. that Laver
(Mathias) forcing adjoins a dominating family of size $\omega_1$.
An immediate consequence of this is that two Laver (Mathias)
reals added iteratively always force $CH$ (Corollary 3.10).
--- Theorem 3.1. enables us to investigate the effect
of adding a Laver real on the invariants in Cicho\'n's diagram
under the assumption that $MA(\sigma$--centered) holds.
Such an assumption (or slightly less) is necessary to ensure
that no cardinals are collapsed. Our results which are
expounded in 3.4. in fact yield an alternative argument for
proving one of the consistency results concerning cardinals
in Cicho\'n's diagram which was originally obtained by
Bartoszy\'nski, Judah and Shelah [BaJS]
using a countable support iteration.
\par
In section 4 we leave our combinatorial considerations for
a while to deal with some descriptive set theory instead.
Recall that a set of reals $A \sub [\omega]^\omega$ is
{\it Ramsey} iff there is $a\in[\omega]^\omega$ so that either
$[a]^\omega \sub A$ or $[a]^\omega \sub [\omega]^\omega
\setminus A$. We shall show that {\it
$\Sigma^1_4$--Mathias--absoluteness} (which means that
$V$ and $V[m]$, where $m$ is Mathias--generic over $V$,
satisfy the same $\Sigma^1_4$--sentences with parameters in
$V$) implies that all $\Sigma^1_3$--sets are Ramsey.
\par
Section 5 studies Miller's rational perfect set forcing
[Mi 2]. It turns out (Theorem 5.1.) that the effect of a Miller
real is rather different from the one of the forcing notions
considered in sections 2 and 3; namely Miller forcing
{\it preserves the axiom $MA(\sigma$--centered}) (which means
that given a model $V$ for "$MA(\sigma$--centered) $+
2^\omega = \kappa$" ($\kappa$ regular uncountable),
"$MA(\sigma$--centered) $+2^\omega = \kappa$" is still
true in $V[m]$, where $m$ is Miller--generic over $V$).
In so far its behaviour is similar to the one of Cohen
forcing (see above). Furthermore it collapses the {\it additivity
of Lebesgue measure} (the smallest size of a family of null
sets the union of which is not null) to $\omega_1$
(Theorem 5.6.).
\bigskip
{\bolds A note for the reader.} We think of Theorems 2.1. and
3.1. (and, to a lesser extent, 5.1.) as the main results of this work.
\par
Sections 2, 3 and 5 are independent of each other, and can be
read without knowledge of the preceding sections (but not
without some sophistication in forcing arguments). 
[However, one
argument in 5.6. is exactly similar to the corresponding 
argument in the proof of 3.4., and is therefore left out.]
The results of section 4 heavily depend on the proof of Theorem 3.1.
and the fact that Mathias forcing satisfies the requirements
of 3.1. (Corollary 3.10), but not on the remainder of
section 3.
\bigskip
{\bolds Notation.} Our notation is fairly standard. We refer
the reader to [Je 1] and [Ku] for set theory in general and
forcing in particular.
\par
Given a finite sequence $\sigma$ (e.g. $\sigma \in \omega^{<\omega}$),
we let $lh(\sigma) := dom (\sigma)$ denote the length of $\sigma$;
for $\ell \in lh(\sigma)$, $\sigma\restrict\ell$ is the restriction
of $\sigma$ to $\ell$. $\hat{\;}$ is used for concatenation
of sequences; and $\la\ra$ is the empty sequence.
Given a tree $T \sub \omega^{<\omega}$, $[T] := \{ f\in\omega^{
\omega} ;\; \forall n \; (f\restrict n \in T ) \}$ denotes
the set of its branches.
\par
A partial order (p.o. for short) is called {\it $\sigma$--centered}
iff $\PP = \bigcup_{n\in\omega} \PP_n$ where each $\PP_n$ is 
{\it centered} (i.e., given $F\sub \PP_n$ finite, there is $q
\in\PP$ so that for all $p\in F \; (q\leq p)$);
and $\PP$ is {\it $\sigma$--linked} iff $\PP = \bigcup_{n\in
\omega} \PP_n$ where each $\PP_n$ is {\it linked}
(i.e., given $p,q \in \PP_n$, there is $r\in \PP$ so that
$r \leq p$ and $r\leq q$). A $\sigma$--centered p.o. is 
$\sigma$--linked; and a $\sigma$--linked p.o. is $ccc$.
An example for a $\sigma$--linked not $\sigma$--centered
p.o. is random forcing, henceforth denoted by $\BB$.
\par
Given a p.o. $\PP \in V$, we shall denote $\PP$--names
by symbols like $\breve f, \breve T, \breve t , ...$ and their
interpretation in $V[G]$ (where $G$ is $\PP$--generic over $V$)
by $\breve f [G], \breve T [G] , \breve t [G] , ...$.
We confuse to some extent Boolean--valued models $V^\PP$ and
forcing extensions $V[G]$, where $G$ is $\PP$--generic over $V$.
For a sentence of the $\PP$--forcing language $\phi$,
$\Bool \phi \Boor$ is the Boolean value of $\phi$.
The symbol $\star$ denotes iteration. --- We say a forcing
notion $\PP$ is {\it generated by a name for a real}
iff there is a $\PP$--name $\breve r$ for an object in
$\omega^\omega$ so that the complete Boolean algebra generated by
the family $\{ \Bool \breve r (i) = n \Boor ;\; i,n \in\omega\}$
equals $r.o.( \PP)$; we express this by saying $\PP =
\PP_{\breve r}$; and we usually denote
this name by the same letter as the 
forcing notion (e.g., $\MM$ for Mathias forcing,
and $\breve m$ for the name of the Mathias--generic real).
\par
Whenever we use $\Sigma^1_n$ or $\Pi^1_n$ we mean the boldface
version.
\bigskip
{\bolds Acknowledgment.} I am very much indebted to Haim Judah
for motivating me to work on the problems considered in this paper
and for various discussions concerning the results.
\Bigskip

{\bold 1. Frontispiece --- forcing notions and cardinal
invariants}
\Smallskip
{\bolds 1.1. Forcing notions.} We shall consider the following
forcing notions. 
\par
\noindent{\it --- Eventually different reals forcing
$\EE$} [Mi 1, section 5]:
\smallskip
\centerline{$(s,G) \in \EE \Longleftrightarrow s \in
\omega^{<\omega} \;\land\; G \in [\omega^\omega]^{<\omega}$}
\smallskip
\centerline{$(s,G) \leq (t,H) \Longleftrightarrow s \supseteq
t \;\land\; G \supseteq H \;\land\; \forall g \in H \; \forall
i \; (dom(t) \leq i < dom (s) \rightarrow s(i) \neq g(i))$}
\smallskip
\noindent{\it --- Localization 
forcing} $\LL$ (see, e.g., [Tr 3, $\S$ 2]):
\smallskip
\centerline{$(\sigma,G) \in \LL \Longleftrightarrow \sigma \in ([
\omega]^{< \omega})^{< \omega} \;\land\; \forall i \in dom(\sigma)
\; (\vert \sigma (i) \vert = i + 1)\;\land\; G \in
[\omega^\omega]^{\leq dom(\sigma) + 1}$}
\smallskip
\centerline{$(\sigma , G) \leq (\tau,H) \Longleftrightarrow
\sigma \supseteq \tau \;\land\; G \supseteq H \;\land\; \forall 
g \in H \;\forall
i \;(dom(\tau) \leq i < dom(\sigma) \rightarrow g(i) \in \sigma
(i))$}
\smallskip
\noindent{\it --- Mathias forcing} $\MM$
[Je 2, part one, section 3]:
\smallskip
\centerline{$(s,S) \in \MM \Longleftrightarrow s \in
\omega^{<\omega}\;\land\;S\in[\omega]^\omega\;\land\;
s$ strictly increasing $\land\; max(ran(s)) < min (S)$}
\smallskip
\centerline{$(s,S) \leq (t,T) \Longleftrightarrow s\supseteq t
\;\land\; S \subseteq T \;\land\; \forall i \in (lh(s)
\setminus lh(t))\;(s(i) \in T)$}
\smallskip
\noindent{\it --- Laver forcing} $\LA$ 
[Je 2, part one, section 3]:
\smallskip
\centerline{$T\in\LA \Longleftrightarrow T\subseteq \omega^{<\omega}
$ is a tree $\land\;\exists\rho\in T\; \forall\sigma\in T \;(
\sigma \subseteq \rho \;\lor\;[\rho\subseteq\sigma \;\land\;
\exists^\infty n\;(\sigma\hat{\;}\langle n \rangle \in T)])$}
\smallskip
\centerline{$T\leq S \Longleftrightarrow T \subseteq S$}
\smallskip
\noindent The $\rho$ required to exist in the above definition
is usually called the {\it stem} of $T$, $stem(T)$. Furthermore,
for $\rho\in T$ we let $succ_T (\rho) := \{ n \in \omega ; \;
\rho\hat{\;}\langle n \rangle \in T \}$, the set of 
{\it successors}
of $\rho$ in $T$, and $T_\rho := \{ \sigma \in T;\; \sigma
\sub \rho \;\lor\; \rho\sub\sigma \}$. We say $T\leq_0 S$
iff $T\leq S$ and $stem (T)=stem(S)$.
\par
Given a non--principal ultrafilter ${\cal U}$ of subsets
of $\omega$, we can define {\it Mathias forcing with
repect to ${\cal U}$}, $\MM_{\cal U}$, as well as {\it
Laver forcing with respect to ${\cal U}$}, $\LA_{\cal U}$:
\sm
\ce{$(s,S) \in \MM_{\cal U} \Longleftrightarrow (s,S)\in
\MM \;\land\; S\in{\cal U}$}
\sm \ce{$T \in \LA_{\cal U} \Longleftrightarrow T\in\LA \;\land\;
\forall\sigma\in T \; (stem (T) \sub \sigma \longrightarrow
succ_T(\sigma) \in{\cal U})$}
\sm
\no The order is the restriction of the Mathias order
(Laver order, respectively).
\par
A non--principal ultrafilter ${\cal U}$ over $\omega$ is called
{\it Ramsey} iff for every partition ${\cal A}$ of $\omega$ so that
$A \not\in {\cal U}$ for all $A\in{\cal A}$, there exists $B\in
{\cal U}$ so that for every $A\in{\cal A}$, $\vert B \cap A \vert
\leq 1$. We note that Mathias forcing $\MM$ is equivalent 
(from the forcing--theoretic point of view) to the two--step
iteration $P(\omega) / fin \star \MM_{\breve{\cal U}}$, where
$(P(\omega) / fin , \leq)$ is the ($\sigma$--closed)
Boolean algebra of subsets of $\omega$ modulo the finite sets
ordered by almost inclusion (i.e. for $A$ and $B \in
[\omega]^\omega$, $A\leq B$ iff $A \sub^* B$ iff $A \setminus
B$ is finite) which generically adjoins a Ramsey ultrafilter
$\breve{\cal U}$ [Ma, 4.9.]. Also, several people have
observed (see [Bl, pp. 238--239] or [JS 1, section 1])
that for Ramsey ultrafilters ${\cal U}$, the forcings $\MM_{\cal U}$
and $\LA_{\cal U}$ are equivalent.
\par
\no{\it --- Miller forcing $\MI$} [Mi 2]:
\sm
\ce{$T\in\MI \Longleftrightarrow T\sub\omega^{<\omega}$ is a 
tree $\land\;\forall\rho\in T\;\exists\sigma\in T \; (\rho
\sub\sigma \;\land\; \exists^\infty n \; (\sigma\hat{\;}\la n
\ra\in T))$}
\sm
\ce{$T\leq S \Longleftrightarrow T\sub S$}
\sm
\no We note that conditions $T$ so that for all $\rho \in T$
either $\exists ! n \; (\rho\hat{\;}\la n \ra \in T)$
or $\exists^\infty n \; (\rho\hat{\;}\la n\ra\in T)$ are dense
in $\MI$, and henceforth restrict our attention to such
conditions. For such a condition $T$ we let $stem(T) =$ the 
unique $\rho\in T$ of minimal length so that $\exists^\infty n
\; (\rho\hat{\;}\la n\ra\in T)$; $split (T) = \{ \rho \in T;
\; \exists^\infty n \; (\rho\hat{\;}\la n \ra \in T)\}$;
and for $\rho\in split(T)$, $succ_T (\rho) = \{ \sigma \in
split (T) ;\; \rho \sub \sigma \;\land$ for no $lh (\rho)
< n < lh(\sigma) \; (\sigma\restrict n \in split (T)) \}$.
\bigskip
{\bolds 1.2. Cardinal invariants.}
Given a $\sigma$-ideal
${\cal I} \subseteq P(2^\omega)$, we let
\smallskip
\itemitem{$add({\cal I}) :=$} the least $\kappa$ such that 
$\exists {\cal F} \in [{\cal I}]^\kappa \; (\bigcup
{\cal F} \not\in {\cal I})$;
\par
\itemitem{$cov ({\cal I}) :=$} the least $\kappa$ such that
$\exists {\cal F} \in [{\cal I}]^\kappa \; (\bigcup
{\cal F} = 2^\omega)$;
\par
\itemitem{$unif({\cal I}) :=$} the least $\kappa$ such that
$[2^\omega]^\kappa \setminus {\cal I} \neq \emptyset$;
\par
\itemitem{$cof ({\cal I}) :=$} the least $\kappa$ such that
$\exists {\cal F} \in [{\cal I}]^\kappa \; \forall A \in {\cal I}
\; \exists B \in {\cal F} \; (A \subseteq B)$.
\smallskip
\noindent We also define
\smallskip
\itemitem{${\bf b} :=$} the least $\kappa$ such that $\exists {\cal F}
\in [\omega^\omega]^\kappa \; \forall f \in \omega^\omega
\; \exists g \in {\cal F} \; \exists^\infty n \;
(g(n) > f(n))$;
\par
\itemitem{${\bf d} :=$} the least $\kappa$ such that $\exists {\cal F}
\in [\omega^\omega]^\kappa \; \forall f \in \omega^\omega \;
\exists g \in {\cal F} \; \forall^\infty n \; (g(n) > f(n))$.
\smallskip
\noindent If ${\cal M}$ is the ideal of meager sets, and
${\cal L}$ is the ideal of null sets, then we can arrange
these cardinals in the following diagram (called Cicho\'n's
diagram).
\bigskip
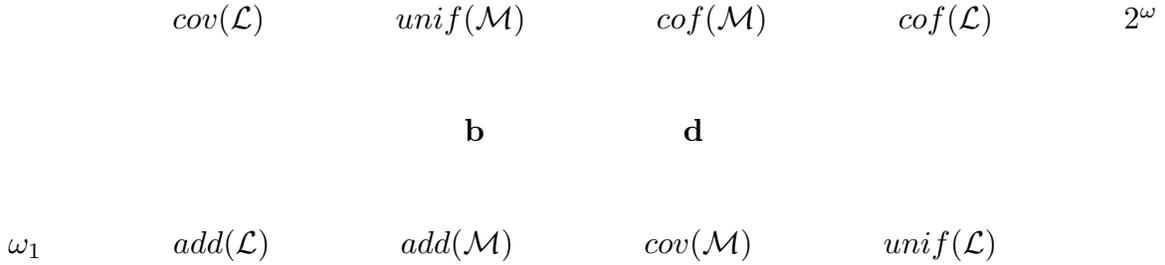
\centerline{{\quad} \Hoskip $cov({\cal L})$ \Hoskip
$unif({\cal M})$ \Hoskip $cof({\cal M})$ \Hoskip $cof({\cal L})$
\Hoskip $2^\omega$}
\Veskip
\centerline{${\bf b}$ \hskip 2.5truecm ${\bf d}$}
\Veskip
\centerline{$\omega_1$ \Hoskip $add({\cal L})$ \Hoskip
$add({\cal M})$ \Hoskip $cov({\cal M})$ \Hoskip $unif({\cal L})$
\Hoskip {\quad}}
\bigskip
\noindent (Here, the invariants grow larger, as one moves up and to the
right in the diagram.) The dotted line says that $add({\cal M})=
min\{ b, cov({\cal M}) \}$ and $cof({\cal M}) = max\{ d,
unif({\cal M}) \}$. For the results which determine the
shape of this diagram, we refer the reader to [Fr]. A survey
on independence proofs showing that no other relations can be
proved between these cardinals can be found in [BaJS].
We shall need the following characterization of the cardinal
$add({\cal L})$, which is due to Bartoszy\'nski [Ba]:
call a function $\phi \in ([\omega]^{<\omega})^\omega$ a 
{\it slalom} iff for all $n\in\omega\;(\vert\phi(n)\vert
\leq n)$; let $\Phi$ be the set of all slaloms; then
\sm
\itemitem{$add({\cal L}) =$} the least $\kappa$ such that
$\exists {\cal F}\in[\omega^\omega]^\kappa \; \forall\phi
\in\Phi\;\exists f\in{\cal F}\;\exists^\infty n \;(f(n)\not\in
\phi(n))$.
\sm
In addition, we shall be interested in the following
cardinals:
\sm
\itemitem{${\bf p}:=$} the least $\kappa$ so that there is
${\cal F}\in [[\omega]^\omega]^\kappa$ with the {\it strong
finite intersection property} (i.e. given finitely many
$A_i \in {\cal F}$, $i<n$, we have $
\vert \bigcap_{i<n} A_i \vert =\omega$) so that $\neg\exists
B\in[\omega]^\omega \;\forall A\in {\cal F}\; (B\sub^* A)$;
\sm
\itemitem{${\bf a}:=$} the least $\kappa$ so that there
is a mad family of size $\kappa$;
\sm
\itemitem{${\bf s}:=$} the least $\kappa$ so that $\exists
{\cal F}\in[[\omega]^\omega]^\kappa \; \forall B
\in[\omega]^\omega \;\exists A \in{\cal F}\; (\vert A \cap B
\vert = \vert B\setminus A\vert =\omega)$;
\sm
\itemitem{${\bf r}:=$} the least $\kappa$ so that $\exists
{\cal F}\in[[\omega]^\omega]^\kappa\;\forall A \in
[\omega]^\omega \;\exists B\in{\cal F} \; (\vert A\cap B\vert
< \omega \;\lor\; \vert B \setminus A \vert < \omega)$;
\sm
\itemitem{${\bf h}:=$} the least $\kappa$ so that there is
${\cal F} \in [[[\omega]^\omega]^{\leq 2^\omega}]^\kappa$
so that all ${\cal A}\in{\cal F}$ are mad and there is no
mad ${\cal B}$ with $\forall {\cal A}\in{\cal F} \;
\forall B \in {\cal B} \; \exists A \in {\cal A} \;
(B \sub^* A )$.
\sm
\par
\no Bell's Theorem says that $2^\omega = {\bf p}$ is equivalent
to $MA(\sigma$--centered) [Be]. It is well--known
that $\omega_1 \leq {\bf p} \leq add({\cal M})$;
that ${\bf p} \leq {\bf h} \leq {\bf s}
\leq unif({\cal M}), unif({\cal L}), {\bf d}$;
that ${\bf h} \leq {\bf b} \leq {\bf a}$; and that ${\bf b}, cov({\cal L}),
cov({\cal M}) \leq {\bf r} \leq 2^\omega$ [vD].
\bigskip
{\bolds 1.3.} The forcings in 1.1. and the cardinals in 1.2. are far from
being unrelated. In fact, some of the former were devised to yield
consistency results concerning the latter using an iterated
forcing construction. However, although it is well--known which
value the cardinals will have after iterating one of the
forcings, it has not yet been investigated what their values are
after adding just one real.
\Bigskip

{\bold 2. The theme --- adding large Luzin sets and small mad families}
\Smallskip
In this section we show that many $ccc$ forcing notions
which have been considered in literature adjoin a Luzin set
of size $2^\omega$ and a mad family of size $\omega_1$.
\par
Assume we have at most countable sets $A_\ell$, $B_\ell$,
natural numbers $n_\ell$, relations $R_\ell
\sub A_\ell \times [B_\ell]^{<\omega}$ and partitions
$A_\ell = \bigcup_{k\in m_\ell} C_{k,\ell}$ for $\ell \in \omega$
such that
\par
\item{(i)} $\vert A_\ell \vert , \vert B_\ell \vert ,
n_\ell , m_\ell$ increase and converge to $\omega$ for
$\ell \to \omega$ (possibly $\vert A_\ell \vert , \vert B_\ell
\vert , m_\ell = \omega$);
\par
\item{(ii)} for all $\ell \in\omega$, if $S\sub T\sub B_\ell
$ are finite and $a\in A_\ell$, then $a R_\ell T$ implies 
$a R_\ell S$;
\par
\item{(iii)} for all $\ell \in \omega$, if $S,T \sub B_\ell
$ are finite, $\vert T \vert \leq n_\ell$ and (necessarily) $S
\not\sub T$, then there is $a \in A_\ell$ with $\neg
(a R_\ell S)$ and $a R_\ell T$;
\par
\item{(iv)} for all $\ell \in \omega$, given $S \sub
B_\ell$ of size $\leq n_\ell$ and $k \in m_\ell$, there is $a \in
C_{k,\ell}$ with $a R_\ell S$.
\par
\no We define a forcing notion $\II$ associated with this
situation as follows.
\sm
\ce{$(\sigma , G) \in \II \Longleftrightarrow \exists \ell
\in\omega\; (\sigma \in \prod_{i\in\ell} A_i \;\land\;
G \in [\prod_{i\in\omega} B_i]^{\leq n_\ell})$}
\sm
\ce{$(\tau,H) \leq (\sigma , G) \Longleftrightarrow
\tau \supseteq \sigma \;\land\; H \supseteq G \;\land\;
\forall i\in dom(\tau ) \setminus dom(\sigma) \;
(\tau(i) R_i \{ g(i) ;\; g\in G \} )$}
\sm
\no We call $\II$ an {\it inflexible} forcing.
Note that $\II$ is $\sigma$--linked (two conditions with
the same initial segment and sets of functions which agree
long enough are compatible by (iv)).
Also both eventually different reals forcing $\EE$ and 
localization forcing $\LL$ are inflexible (for both,
let $B_\ell = \omega$, $\omega = \bigcup_{k\in\omega}
C_k$ any partition into infinitely many infinite pieces,
and $n_\ell = \ell$; for $\EE$, $A_\ell = \omega$,
$C_{k,\ell} =C_k$, and $R_\ell = \not\in$;
for $\LL$, $A_\ell = [\omega]^{\ell + 1}$,
$C_{k,\ell} = \{ a \in A_\ell ;\; max (a) \in
C_k \}$, and $R_\ell = \supseteq$).
\sm
{\bolds 2.1. Theorem.} {\it An inflexible forcing $\II$ adds
a Luzin set of size $2^\omega$.}
\sm
{\capit Proof.} Let $\{ f_{\alpha , m ,i} ;\; \alpha \in 
2^\omega \;\land\; m \in \omega \;\land\; i \in n_m \}
\sub \prod_{i\in\omega} B_i$ be a family of almost disjoint
functions (i.e. $(\alpha , \ell , i) \neq (\beta , m ,j)$
implies $\forall^\infty n \; (f_{\alpha , \ell , i}  (n)
\neq f_{\beta , m, j} (n))$ --- such a family is constructed
like an a.d. family of size $2^\omega$).
We shall define a sequence $\la \breve f_\alpha ;\;
\alpha\in 2^\omega \ra$ of $\II$--names for objects in
$\omega^\omega$.
\par
Fix $\alpha \in 2^\omega$ and $m \in \omega$. Note that
\sm
\ce{$D_{\alpha,m} := \{ (\sigma, \{f_{\alpha,m,i};\; i\in n_m \}
); \; \exists \ell \geq m \; (\sigma \in \prod_{i\in\ell}
A_i \;\land$}
\par \ce{$\land\; (\ell > m \longrightarrow \neg (\sigma (\ell
- 1) R_{\ell - 1} \{ f_{\alpha ,m,i} (\ell - 1) ;
i \in n_m \}))) \}$}
\sm
\no is a maximal antichain in $\II$. For $(\sigma ,
\{ f_{\alpha , m ,i} ; \; i \in n_m \}) \in
D_{\alpha, m}$, we set
$$ (\sigma , \{ f_{\alpha , m,i} ;\; i \in n_m \})
\forces_\II "\breve f_\alpha (m) = \ell"
\longleftrightarrow \cases{lh(\sigma) \geq 2 \; {\rm
and } \; \sigma (lh(\sigma) - 2) \in C_{\ell ,
lh (\sigma) - 2} & or \cr
lh(\sigma) < 2 \; {\rm and} \; \ell =0. &\cr}$$
This completes the definition of the sequence
of names $\la \breve f_\alpha ; \; \alpha < 2^\omega
\ra$.
\par
Next let $\breve T$ be an $\II$--name so that
\sm
\ce{$\forces_\II " \breve T \sub \omega^{<\omega}$ is a
nowhere dense tree".}
\sm
\no
This means that given $s \in \omega^{<\omega}$
we have an $\II$--name $\breve t_s$ so that
\sm
\ce{$\forces_\II "s \sub \breve t_s$ and $\breve t_s
\not\in \breve T"$.}
\sm
\no Let $I_s := \{ (\sigma^n_s , G^n_s ) ;\; n\in\omega\}$
be a maximal antichain of conditions deciding $\breve t_s$.
We say $\alpha \in 2^\omega$ is {\it $(m,s,n)$--bad}
iff
\sm
\ce{$\forall^\infty i\in\omega\; (\vert \{ f_{\alpha ,m,j} (i) ;
\; j \in n_m \} \setminus \{ g(i) ; \; g \in G^n_s \}
\vert \leq {n_m \over 2})$.}
\sm
\no By almost--disjointness of the functions $f_{\alpha ,m,j}$
only finitely many $\alpha $ can be $(m,s,n)$--bad. Let $
{\cal A}:= \{ \alpha ;\; \exists (m,s,n) \in \omega\times
\omega^{<\omega}\times\omega \; (\alpha$ is $(m,s,n)$--
bad$)\}$. ${\cal A}$ is an at most countable
subset of $2^\omega$. To finish the proof of Theorem 2.1. it
suffices to show:
\sm
{\bolds 2.2. Claim.} {\it For all $\alpha \in 2^\omega \setminus
{\cal A} $, $\forces_\II "\breve f_\alpha \not\in \breve T"$.}
\sm
{\capit Proof of Claim.} Suppose, by way of contradiction, that
there are $\alpha \in 2^\omega \setminus {\cal A}$ and $(\sigma
, G) \in \II$ so that
\sm
\ce{$(\sigma , G) \forces_\II " \breve f_\alpha \in \breve T".$}
\sm
\no Choose $m\in\omega$ so that $\vert G \vert < {n_m \over 2}$.
Let $G' := G \cup \{ f_{\alpha , \ell ,i} ;\; \ell < m \;\land\;
i \in n_\ell \}$. Extending $\sigma$, if necessary, we may assume
that $(\sigma , G')$ decides $\breve f_\alpha \restrict m$,
say $(\sigma , G') \forces_\II "\breve f_\alpha \restrict m
= s"$. For some $n\in\omega$, $(\sigma , G')$ and $(\sigma^n_s ,
G^n_s)$ are compatible with common extension $(\tau, G'')$;
without loss $G''$ is of the form $G'' = G' \cup G^n_s$
(cf (ii) in the definition of inflexibility). Let
$t \in \omega^{<\omega}$ be such that $(\sigma^n_s, G^n_s)
\forces_\II "\breve t_s = t"$. Set $m': = lh(t)$;
without loss $m' > m$ (otherwise we are done).
\par
For $m' > \ell \geq m$ we recursively choose $i_\ell$
so that
\par
\item{1)} $i_\ell \geq i_{\ell -1} + 2$ (where $i_{m-1} + 1
= lh (\tau)$) and $i_\ell \geq \ell$;
\par
\item{2)} $\vert \{ f_{\alpha , \ell , j} (i_\ell) ; \;
j \in n_\ell \} \setminus \{ g(i_\ell) ; \; g \in
G^n_s \} \vert > {n_\ell \over 2}$;
\par
\item{3)} $\{ f_{\alpha ,\ell , j} (i_\ell) ;\; j \in n_\ell
\} \cap \{ f_{\alpha , \ell ' , j} (i_\ell) ;\;
\ell ' < \ell \;\land\; j\in n_{\ell '} \} = \emptyset$;
\par
\item{4)} $t(\ell) < m_{i_\ell - 1}$;
\par
\item{5)} $n_{i_\ell + 1} \geq \vert G'' \vert + \sum_{\ell
\geq \ell ' \geq m} n_{\ell '}$.
\par
\no Note that 2) is possible because $\alpha \not\in{\cal A}$
(the other conditions can be trivially satisfied by choosing
$i_\ell$ large enough). 
Next, again by recursion on $m' > \ell \geq m$ we extend $\tau$
to $\tau_\ell$ of length $i_\ell + 1$ so that
\par
\item{a)} $\tau \sub \tau_{\ell - 1} \sub \tau_\ell$
(where $\tau_{m-1} = \tau$);
\par
\item{b)} $\forall i \; (i_{\ell - 1} + 1 \leq i \leq i_\ell 
\longrightarrow (\tau_\ell (i) R_i (\{ g(i) ;\; g \in G''\}
\cup \{ f_{\alpha , \ell ' , j} (i) ;\; \ell ' < \ell \;\land\;
j \in n_{\ell '}\}) ))$;
\par
\item{c)} $\tau_\ell (i_\ell - 1) \in C_{t(\ell) , i_\ell - 1}$;
\par
\item{d)} $\neg (\tau_\ell (i_\ell ) R_{i_\ell} \{ f_{\alpha ,
\ell,
j} (i_\ell) ;\; j \in n_\ell \} )$.
\par
\no This is possible by (iii) and (iv) in the definition of
inflexibility. But $(\tau_{m'-1} , G'' \cup \{ f_{\alpha , \ell , j}
;\; \ell < m' \;\land\; j \in n_\ell \} ) \leq (\tau , G'')
\leq (\sigma , G)$ and
\sm
\ce{$(\tau_{m' - 1} , G'' \cup \{ f_{\alpha , \ell ,j}
;\; \ell < m' \;\land\; j \in n_\ell \} ) \forces_\II "
\breve f_\alpha \restrict m' = t = \breve t_s \not\in
\breve T"$,}
\sm
\no a contradiction. $\qed$ $\qed$
\bigskip
Since the following result is proved by a combination of the methods
of the proof of 2.1. and the one of [BJS, Theorem 2.2.],
we merely sketch the argument.
\sm
{\bolds 2.3. Theorem.} {\it An inflexible forcing $\II$ adjoins
a mad family of size $\omega_1$.}
\sm
{\capit Proof.} We start with an easy observation from 
[BJS, 2.2.].
\sm
{\bolds 2.4. Observation.}
{\it Let $\langle N_\alpha ; \; \omega
\leq \alpha < \omega_1 \rangle$, $\langle h_\alpha ; \; \omega
\leq \alpha < \omega_1 \rangle$ and $\langle r_\alpha ; \;
\omega \leq \alpha < \omega_1 \rangle$ be sequences such
that $N_\alpha \prec H(\chi)$ is countable and $N_\alpha
\prec N_\beta$ for $\alpha < \beta$, $h_\alpha \in \alpha^\omega
\cap N_\alpha$ is one-to-one and onto, $r_\alpha \in \omega^\omega$
is Cohen over $N_\alpha$ and $\langle r_\alpha ; \; \alpha <
\beta \rangle \in N_\beta$. Define recursively sets $C_\alpha$
for $\alpha < \omega_1$. $\langle C_n ; \; n \in \omega \rangle
$ is a partition of $\omega$ into countable pieces lying in $N_\omega$.
For $\alpha \geq\omega$, $C_\alpha := \{ r_\alpha (n) ; \; n \in
\omega \; \land \; \forall m < n \; (r_\alpha (n) \not\in C_{h_\alpha
(m)} ) \}$. Then $\{ C_\alpha ; \; \alpha \in \omega_1 \}$
is an a. d. family.} $\qed$
\smallskip
As in the proof of 2.1., we let $\{ f_{\alpha , m,i} ;\;
\alpha \in\omega_1 \;\land\; m \in \omega \;\land\; i\in n_m \}
\sub \prod_{i\in\omega} B_i$ be a family of almost disjoint
functions; and we define the sequence $\la \breve f_\alpha ; \;
\alpha \in \omega_1 \ra$ of $\II$--names for functions from
$\omega$ to $\omega$ as before. By Theorem 2.1., $\la\breve f_\alpha
[i] ;\; \alpha \in \omega_1 \ra$ is Luzin in $V[i]$ (where $i$ is 
the $\II$--generic real), and thus we can find (still
in $V[i]$) a strictly increasing
function $\phi : \omega_1 \setminus \omega \to \omega_1$ and
sequences $\la N_\alpha ;\; \omega \leq \alpha < \omega_1 \ra$,
$\la h_\alpha ;\; \omega \leq \alpha < \omega_1 \ra$ so that
for $r_\alpha := \breve f_{\phi (\alpha) } [i]$ the requirements
of 2.4. are satisfied. By $ccc$--ness of $\II$, we may assume
that $\phi \in V$; and, in fact, $\phi = id$. Also note
that we can suppose that $\la h_\alpha ; \; \omega \leq
\alpha < \omega_1 \ra\in V$. We shall prove that
the resulting family $\la C_\alpha ;\; \alpha < \omega_1 \ra$
is indeed mad.
\par
For suppose not, and let $\breve C$ be an $\II$--name so that
\sm
\ce{$\forces_\II "\forall \alpha < \omega_1 \;
(\vert \breve C_\alpha \cap \breve C \vert < \omega)".$}
\sm
\no Let $\breve f$ be an $\II$--name for the strictly increasing
enumeration of $\breve C$. Similar to the proof of Theorem 2.1., let
$I_m := \{ (\sigma^n_m , G^n_m ) ;\; n
\in\omega \}$ be a maximal antichain of conditions deciding
$\breve f (m)$. Next, for $\alpha < \omega_1$, let $\breve k_\alpha$
be an $\II$--name for a natural number so that
\sm
\ce{$\forces_\II " \breve C_\alpha \cap \breve C \sub \breve
k_\alpha "$.}
\sm
\no Let $I_\alpha ' := \{ (\sigma^{\prime n}_\alpha , 
G^{\prime n}_\alpha  )
;\; n\in\omega\}$ be a maximal antichain of conditions deciding
$\breve k_\alpha$. For $\Gamma \subseteq \omega_1$ finite and
$m, m' , N \in\omega$ say $\alpha$ is {\it $(m,m',\Gamma, N)$--bad
} iff 
\sm
\ce{$\forall^\infty i \in\omega \; (\vert \{ f_{\alpha , m,j}
(i) ;\; j\in n_m \} \setminus \{ g(i) ; \; \exists n < N \;
(\exists \beta \in \Gamma \; (g \in G_\beta^{\prime n} ) \;\lor\;
g \in G^n_{m'} ) \} \vert \leq {n_m \over 2}).$}
\sm
\no Only finitely many $\alpha$ can be $(m,m',\Gamma ,N)$--bad.
Finally we choose $\alpha < \omega_1$ so that \par
\item{$(\clubsuit)$} if $\Gamma \sub \alpha$ and $m, m' , N \in
\omega$ and $\beta $ is $(m, m' ,\Gamma, N)$--bad, then
$\beta < \alpha$.
\par
\no Choose $n \in \omega$ arbitrarily, and let $k\in\omega$ be such
that $(\sigma^{\prime n}_\alpha  , G^{\prime n}_\alpha ) 
\forces_\II "\breve k_\alpha
= k"$. We reach a contradiction by showing:
\sm
{\bolds 2.5. Claim.} {\it $(\sigma^{\prime n}_\alpha , G^{\prime n}_\alpha )
\not\forces_\II "\breve C_\alpha \cap \breve C \sub k"$.}
\sm
{\capit Proof of Claim.} Choose $m\in\omega$ so that $\vert
G^{\prime n}_\alpha \vert < {n_m \over 2}$. Let $G' := 
G^{\prime n}_\alpha
\cup \{ f_{\alpha ,\ell, i } ;\; \ell < m \;\land\; i \in 
n_\ell \}$. Extend $\sigma^{\prime n}_\alpha$ to $\sigma$, if necessary,
so that $(\sigma, G') \leq (\sigma^{\prime n}_\alpha ,
G^{\prime n}_\alpha)$ is a condition and decides $\breve f_\alpha
\restrict m$, say $(\sigma , G') \forces_\II "\breve f_\alpha 
\restrict m = s$". Let $\Gamma := \{ h_\alpha (\ell) ; \;
\ell < m \}$; for each $\beta \in \Gamma$ find (recursively)
$n_\beta \in\omega$ so that $\{ (\sigma^{\prime n_\beta}_\beta,
G^{\prime n_\beta}_\beta); \; \beta \in\Gamma \} \cup \{(\sigma ,
G')\}$ has a common extension, say $(\tau, G'')$; let $k_\beta
\in\omega$ be so that $(\sigma^{\prime n_\beta}_\beta, G^{\prime
n_\beta}_\beta
) \forces_\II "\breve k_\beta = k_\beta"$. Let $m'$ be larger than
the maximum of the $k_\beta$ and $k$. Find $n' \in\omega$ so that
$(\tau , G'')$ and $(\sigma^{n'}_{m'} , G^{n'}_{m'} )$ have a
common extension, say $(\tau ' ,G''')$; without loss $G'''$
is of the form $G''' = G' \cup \bigcup_{\beta \in \Gamma}
G_\beta^{\prime n_\beta} \cup G^{n'}_{m '}$; let $\tilde\ell\in\omega$ be
such that $(\sigma^{n'}_{m'} , G^{n'}_{m'}) \forces_\II
"\breve f (m') = \tilde \ell "$.
As $\breve f$ is an $\II$--name for a strictly increasing
function, $\tilde\ell \geq m' \geq k$.
\par
Using an argument similar to the final argument in the proof of
Claim 2.2. (applying that $\alpha$ is not $(m,m',\Gamma, N)$--bad
where $N$ is larger than the maximum of the $n_\beta$ and $n'$
$(\clubsuit)$)
we can find $(\tau '' , G ''' \cup \{ f_{\alpha , m , i} ;\;
i \in n_m \} ) \leq (\tau ' , G''') \leq (\sigma^{\prime n}_\alpha 
, G^{\prime n}_\alpha)$ so that
\sm
\ce{$(\tau '' , G''' \cup \{f_{\alpha ,m,i} ;\; i \in n_m \} )
\forces_\II " \breve f_\alpha (m) =\tilde \ell = \breve f(m')
\in (\breve C \cap \breve C_\alpha ) \setminus k "$,}
\sm
\no a contradiction (note that $\breve f_\alpha (m)$ is forced
to belong to $\breve C_\alpha$, because a weaker condition forces
it not to belong to any $\breve C_{h_\alpha (\ell)}$,
$\ell < m$, see 2.4.). $\qed$ $\qed$
\bigskip
As the corresponding result in section 2 of [BJS], Theorem
2.1. has consequences concerning the values of the invariants
in Cicho\'n's diagram after adding a real via an inflexible
forcing notion.
\sm
{\bolds 2.6. Corollary.} {\it (a) Let $\II \in V$ be inflexible;
then $V^\II \models  " unif ({\cal M}) = \omega_1 $ and $cov
({\cal M}) = 2^\omega"$. Thus the invariants on the left--hand
sice in Cicho\'n's diagram all equal $\omega_1$, whereas
those on the right--hand side are equal to $2^\omega$ in 
$V[i]$. Furthermore, ${\bf p} = {\bf h} = {\bf s} = {\bf a}
= \omega_1$ and ${\bf r} = 2^\omega$ in $V[i]$.
\par
(b) In particular, after adding an eventually different real $e$
or a localization real $\ell$, $unif({\cal M}) = \omega_1$ and
$cov({\cal M}) = 2^\omega$ in $V[e]$ (or $V[\ell]$,
respectively). $\qed$}
\sm
\no Part (b) improves earlier partial results in 
[Br 2, $\S$ 1]. \par
We note that our definition of {\it inflexible} is general
enough to provide us with the same results for
the forcings adding an eventually different or a localization
real below a given function $g\in\omega^\omega$ converging
to infinity (in case of localization, this forcing has been
considered by Pawlikowski, see [Pa 1, section 2] for details).
\bigskip
We close this section with some comments concerning the
structure of inflexible forcing notions. We say a p.o.
$\PP$ has {\it $(\kappa , \omega)$--caliber} iff for every $A
\sub \PP$ of size $\kappa$ there is a countable $B \sub A$
and a $p\in\PP$ with $\forall q \in B \; (p \leq q)$. This
notion is implicit in [Tr 2] (see also [Br 2, section 1]).
Using an argument with almost disjoint functions as before we
see that an inflexible forcing $\II$ does not have $(2^\omega ,
\omega)$--caliber. This suggests that we ask:
\sm
{\bolds 2.7. Question.} {\it Assume $\PP$ is a $ccc$ forcing
notion generated by a name for a real of size $\kappa$,
$\kappa \leq 2^\omega$, which does not have $(\kappa ,
\omega)$--caliber. Does this imply that $\PP$ adjoins a Luzin
set of size $\kappa$?}
\sm
\no We note that we cannot drop any of our hypotheses. We
shall see non--$ccc$ counterexamples in the next section;
the algebra $\BB_\kappa$ adding $\kappa$ random reals is a
counterexample (it adds a Sierpi\'nski set [i.e. an uncountable
set of reals which has at most countable intersection with every null
set], but not
a Luzin set, of size $\kappa$), which is not generated by a
name for a real; finally the random algebra $\BB$ itself is a counterexample 
having $(2^\omega, \omega)$--caliber.
\Bigskip
%

{\bold 3. Variation --- adjoining small dominating families}
\Smallskip
This section is dedicated to the investigation of several proper
forcing notions which add a dominating real like Mathias and
Laver forcing.
\par
We say a forcing notion $\PP$ is {\it Laver--like} iff $\PP
=\PP_{\breve r}$ where $\breve r$ is a $\PP$--name for a dominating
real and given $p\in\PP$ there is a Laver tree
$T\sub\omega^{<\omega}$ so that $\forall \sigma \in T \;
(p(T_\sigma) := \bigcap_{n\in\omega} \bigcup_{\tau\in T_\sigma ,
lh(\tau) = n} (p \cap \Bool \breve r \restrict lh (\tau) =
\tau \Boor  ) \in r.o.(\PP) \setminus \{ 0 \} )$
[we usually express this in a somewhat sloppy way by saying
$p(T) \neq 0$ where $p(T) = p(T_{stem(T)})$].
Note that all forcings adding a dominating real 
which we have encountered so far have this property. Following
Judah and Shelah ([JS 1, 1.5.] and [JS 2, 1.9.]), for a Laver
tree $T$, we say $A\sub T$ is a {\it front} iff $\sigma \neq \tau$
in $A$ implies $\sigma \not\sub \tau$ and for all $f \in[T]$
there is $n\in\omega$ so that $f\restrict n \in A$.
We say a Laver--like $\PP$ has {\it strong fusion} iff given 
countably many open dense sets $D_n \sub \PP$ and $p\in\PP$,
there is a Laver tree $T$ so that $p(T) \neq 0$ and for each $n$,
$\{ \sigma \in T ;\; p(T) \cap \Bool \breve r \restrict
lh(\sigma) = \sigma \Boor \in D_n \}$ contains a front. Finally,
a Laver--like $\PP$ is {\it closed under finite changes}
iff given $p\in\PP$, Laver trees $T$ and $T' \sub T$ so that 
for all $\sigma \in T' \; (\vert succ_T (\sigma) \setminus
succ_{T'} (\sigma) \vert < \omega )$, if $p(T) \neq 0$,
then $p(T') \neq 0$. We call $\FF$ a {\it flexible}
forcing notion iff $\FF$ is Laver--like, has strong fusion
and is closed under finite changes.
\sm
{\bolds 3.1. Theorem.} {\it A flexible forcing $\FF$ adds a
dominating family of size $\omega_1$.}
\sm
{\capit Proof.} 
We start by creating names for dominating functions.
Fix $\alpha \in\omega_1$. We define the tree $T_\alpha
\sub\omega^{<\omega}$ and the function $\rho_\alpha :
T_\alpha \to \alpha +1 $ by recursion on the levels.
\par
\item{$\odot$} $\la\ra\in T_\alpha$, $\rho_\alpha
(\la\ra)=\alpha$.\par
\item{$\odot$} Assume $\sigma \in T_\alpha$, $\alpha \geq
\rho_\alpha (\sigma) > 0$. Then $\forall n \; (\sigma \hat{\;}
\la n\ra \in T_\alpha )$. If $\rho_\alpha (\sigma) = \beta + 1$
is a successor, $\forall n \; (\rho_\alpha (\sigma \hat{\;} \la n\ra ) =
\beta)$. If $\rho_\alpha (\sigma)$ is a limit choose a strictly
increasing sequence $\la \beta_n ;\; n\in\omega\ra$ such that
$\rho_\alpha (\sigma) = \bigcup_{n\in\omega} \beta_n$, and
let $\forall n \; (\rho_\alpha (\sigma \hat{\;} \la n \ra ) =
\beta_n)$.
\par
\item{$\odot$} Assume $\sigma \in T_\alpha$, $\rho_\alpha 
(\sigma) = 0$. Then $\forall n \; (\sigma \hat{\;} \la n \ra
\not\in T_\alpha )$.
\par
\no This concludes the definition of $T_\alpha $ and $\rho_\alpha$.
It is immediate that $T_\alpha$ has no infinite branches, and that
$\rho_\alpha$ is the canonical rank function witnessing this.
\par
Next we recursively introduce a sequence $\la \rho^n_\alpha
;\; n\in\omega\ra$ of functions defined on subsets of
$\omega^{<\omega}$:
\par
\item{$\bullet$} $dom \rho^0_\alpha = \omega^{<\omega}$,
$$\rho^0_\alpha (\sigma) =\cases{\rho_\alpha (\sigma) 
& if $\sigma \in T_\alpha$ \cr 0 & otherwise \cr}$$
\par
\item{$\bullet$} $dom \rho^n_\alpha = \{ \sigma ; \;
\rho^{n-1}_\alpha (\sigma) = 0 \}$, for $\sigma \in
dom \rho^n_\alpha$ choose $m\in\omega$ minimal 
such that $\rho^{n-1}_\alpha (\sigma \restrict m ) =0$
and assume $\sigma = \sigma \restrict m \hat{\;}
\tau$. Then
$$\rho_\alpha^n (\sigma) = \cases{ \rho_\alpha (\tau)
& if $\tau \in T_\alpha$ \cr 0 & otherwise}$$
(Note $\rho^n_\alpha (\sigma) = \alpha \longleftrightarrow
\tau =\la\ra \longleftrightarrow m = lh(\sigma)$.)
\par
\no This concludes the definition of $\la \rho^n_\alpha ;\; 
n\in\omega\ra$. Now it is easy to define an $\FF$--name 
$\breve f_\alpha$ for a function in $\omega^\omega$ 
as follows: $p\in\FF$ decides the value $\breve f_\alpha
(n)$ iff there is $\sigma \in \omega^{<\omega}$ so that
$p \leq \Bool \breve r \restrict lh (\sigma) = \sigma \Boor$
and $\rho^n_\alpha (\sigma) = 0$; in this case let $m\in
\omega$ be minimal so that $\rho^n_\alpha (\sigma
\restrict m ) = 0$; then we say
\sm
\ce{$p\forces_\FF " \breve f_\alpha (n) = \sigma (m-1)".$}
\sm
{\bolds 3.2. Claim.} {\it $\forces_\FF " \{ \breve f_\alpha
;\; \alpha < \omega_1 \}$ is a dominating family$"$.}
\sm
{\capit Proof of Claim.} Let $\breve g$ be an $\FF$--name for a
function from $\omega$ to $\omega$. Let $D_n = \{ p\in\FF ;\; p$
decides $\breve g (n) \}$. Given $p \in \FF$, apply strong
fusion to get a Laver tree $T$ so that $p\geq p(T) \neq 0$
and for each $n$, $\{ \sigma \in T ;\; p(T) \cap \Bool
\breve r \restrict lh(\sigma) = \sigma \Boor \in D_n \}$
contains a front $(*)$. For $n\in\omega$ we define a
function $rk_n : T \to \omega_1$ by recursion on the 
ordinals as follows:
\sm
\ce{$rk_n(\eta) = 0 \longleftrightarrow p(T) \cap \Bool \breve r
\restrict lh(\eta) = \eta \Boor \in D_n$}
\sm
\ce{$rk_n (\eta) = sup \{ rk_n (\eta \hat{\;} \la m\ra )
+ 1 ; \; \eta \hat{\;} \la m\ra \in T \}$ otherwise}
\sm
\no It follows from $(*)$ above that $rk_n$ is defined on
all of $T$. Let $\tilde \alpha := sup \{ rk_n (\la\ra)
+1 ;\; n\in\omega\}$.
\sm
{\bolds 3.3. Subclaim.} {\it $p(T) \forces_\FF "\forall \alpha
\geq \tilde\alpha \;\forall^\infty n \; (\breve f_\alpha (n) >
\breve g(n))"$.}
\sm
{\capit Proof of Subclaim.} Choose $\alpha \geq \tilde\alpha$
and a condition $q\leq p(T)$; without loss
$q$ is of the form $q(S)$ for some Laver tree $S\sub T$.
Let $\sigma = stem(S)$. Let $n$ be minimal so that either
$\rho_\alpha^n (\sigma) = \alpha$ or $\sigma \not\in
dom \rho^n_\alpha$. Extending $\sigma$ if necessary, we 
can assume that the former case holds. We construct by recursion 
on the levels $S' \sub S$ with $stem (S') = \sigma$ so that
for all $\tau \in S' \; (\vert succ_S (\tau) \setminus
succ_{S'} (\tau) \vert < \omega )$ and
\sm
\ce{$q(S') \forces_\FF "\forall m \geq n \; (\breve f_\alpha
(m) > \breve g(m))"$. $(**)$}
\sm
\no
(Note that $q(S') \neq 0$ by closure under finite changes.)
We guarantee along the construction that if $\tau\in S'$
then for all $m\geq n$ either $0 \leq rk_m (\tau) < \rho^m_\alpha
(\tau)$ or $(0 = rk_m (\tau) = \rho^m_\alpha (\tau)$ and
some $\ell_1 < \ell_2 \leq lh(\tau)$ are minimal with
$rk_m (\tau \restrict \ell_1) = 0$ and $\rho^m_\alpha (\tau \restrict
\ell_2) = 0$ and $\tau(\ell_2 - 1) > \ell_3$ where $p(T)
\cap \Bool \breve r \restrict \ell_1 = \tau \restrict \ell_1
\Boor \forces_\FF " \breve g (m) = \ell_3")$.
(this clearly implies $(**)$).
\par
To do the construction assume we have put $\tau\supseteq \sigma$,
$\tau \in S$, into $S'$. We do not have to care about $m$'s
so that $\rho^m_\alpha (\tau) = 0$. Assume, therefore, that $m \geq
n$ is minimal so that $\rho^m_\alpha (\tau) > 0$; then
also $\rho^m_\alpha (\tau) > rk_m (\tau)$. In case $\rho^m_\alpha
(\tau) = 1 > 0=rk_m (\tau) $ find $\ell_1 \leq lh(\tau)$
minimal such that $rk_m (\tau \restrict \ell_1) = 0$ and
$\ell_3$ such that $p(T) \cap \Bool \breve r \restrict \ell_1
= \tau \restrict \ell_1 \Boor \forces_\FF " \breve g (m) =
\ell_3"$. Put $\tau \hat{\;} \la \ell\ra$ into $S'$ iff
$\tau\hat{\;}\la\ell\ra\in S$ and $\ell > \ell_3$. In all
other cases $\rho^m_\alpha (\tau \hat{\;} \la\ell\ra ) >
rk_m (\tau\hat{\;}\la\ell\ra )$ for almost all $\ell$
such that $\tau\hat{\;}\la\ell\ra\in S$; and we can put these
into $S'$. Finally, we do not have to care about $m' > m$
because $\tau\not\in dom \rho^{m'}_\alpha$. This concludes
the construction, the proof of the Subclaim, Claim, and
Theorem. $\qed $ $\qed$ $ \qed$
\bigskip
We can now determine the effect of adding a Laver
real on almost all invariants in Cicho\'n's diagram under $MA(
\sigma$--centered).
\sm
{\bolds 3.4. Theorem.} {\it Let $V \models ZFC$. Then \par
(a) $V^\LA \models " {\bf d} = \omega_1 "$. \par
(b) If $V \models " MA(\sigma$--centered$) + 2^\omega =
\kappa "$, then $V^\LA \models " unif ({\cal L}) = 
unif ({\cal M}) = \kappa = 2^\omega "$. \par
(c) Assume $V\models GCH$, let $\kappa$ be regular uncountable,
and assume $\la \PP_\alpha , \breve \QQ_\alpha ;\; \alpha < \kappa
\ra$ is an iterated forcing construction with finite
supports so that
\par
\ce{$\forces_{\PP_\alpha} " \breve \QQ_\alpha$ is $\sigma$--centered
$"$ and}
\par
\ce{$\forces_{\PP_\kappa} "MA(\sigma$--centered)$"$.}
\no Then $V^{\PP_\kappa \star \breve\LA} \models "unif({\cal
L}) = unif ({\cal M}) = \kappa = 2^\omega + {\bf d}
= cov({\cal L}) = \omega_1 "$.}
\sm
{\capit Proof.} (a) To show ${\bf d} = \omega_1$, we will
apply the previous Theorem. Laver forcing clearly is Laver--like
and closed under finite changes, and strong fusion was
proved by Judah and Shelah in [JS 2, 1.9.]. We repeat the
argument here for completeness' sake.
\par
Take a Laver tree $T$ and $D\sub\LA$ open dense. For each 
$\sigma \in T$ we define the {\it rank} of $\sigma$,
$\rho (\sigma)$, by recursion on the ordinals:
\sm
\ce{$\rho(\sigma) = 0 \longleftrightarrow \exists S
\leq_0 T_\sigma \; (S \in D)$.}
\sm
\ce{$\rho(\sigma) = \alpha \longleftrightarrow $ for no
$\beta < \alpha \; (\rho (\sigma) = \beta)$ and $\exists^\infty
n \; (\sigma \hat{\;} \la n \ra \in T$ and $\rho(\sigma \hat{\;}
\la n\ra ) < \alpha)$.}
\sm
\no If $\rho(\sigma)$ is undefined, we put $\rho (\sigma) =
\infty$.
\par
We claim that $\rho (stem (T)) < \infty$. Otherwise build
a condition $T' \leq T$ recursively so that if $\sigma \in T'$
then $\rho (\sigma) = \infty$: $stem (T) \in T'$; by definition
of $\rho$, almost all successors of $stem(T)$ in $T$ have
rank $\infty$; put them into $T'$ and proceed. In the end choose
$S' \leq T'$ so that $S' \in D$. Then $\rho (stem (S'))
=0$ by definition; but also $\rho(stem(S')) = \infty$ by
construction of $T'$, a contradiction.
\par
Next, build $T' \leq_0 T$ recursively so that if $\sigma \in T'$
and $\rho (\sigma) > 0$, then for all $n\in\omega$ with $\sigma
\hat{\;} \la n\ra\in T'$, $\rho (\sigma \hat{\;} \la n\ra ) <
\rho (\sigma)$, and so that if $\sigma \in T'$, $\rho (\sigma) =
0$ and for all $n < lh(\sigma)$ $\rho (\sigma \restrict n ) >
0$, then $T'_\sigma \in D$.
\par
This gives a front for one open dense set. For countably many just
use one of the standard fusion arguments for Laver forcing.
\sm
(b) For the uniformities we have to modify the argument
that $\LA$ doesn't collapse cardinals under $MA(\sigma$--centered),
due to Judah, Miller and Shelah [JMS, section 3].
\par
Following [JMS, section 3], we let $Q = \{ \bar A = \la
A_\sigma ;\; \sigma\in\omega^{<\omega}\ra ;\; \forall\sigma\;
(A_\sigma \in [\omega]^\omega ) \}$. With $\sigma \in \omega^{<\omega}
$ and $\bar A = \la A_\sigma ;\; \sigma\in\omega^{<\omega} \ra\in
Q$ we can associate a condition $T_\sigma = T_\sigma(\bar A)
\in \LA$ as follows. $T_\sigma (\bar A)$ is the unique Laver tree
with $stem (T_\sigma) =\sigma$ and if $\tau\supseteq\sigma$
and $\tau\in T_\sigma$ then $succ_{T_\sigma} (\tau) = A_\tau$.
We put $\bar A = \la A_\sigma ;\;\sigma\in\omega^{<\omega}\ra
\leq \bar B =\la B_\sigma ;\; \sigma \in\omega^{<\omega} \ra$ iff
$\forall\sigma \; (A_\sigma \sub B_\sigma)$ iff $\forall \sigma
\; (T_\sigma (\bar A) \leq T_\sigma (\bar B))$.
\par
Let $\la I_n ;\; n\in\omega\ra$ be a partition of $\omega$
into finite intervals so that $max(I_n) + 1 = min (I_{n+1})$
and $\vert I_n \vert \geq n^2$. The following is proved easily
using the Laver property and the argument of [JMS, Lemma 3.1.].
\sm
{\bolds 3.5. Lemma.} {\it Assume $\forces_\LA " \breve f \in
2^\omega "$ and $\bar B = \la B_\sigma ;\; \sigma
\in\omega^{<\omega} \ra \in Q$. Then there are $\phi_\sigma \in 
\prod_n [2^{I_n}]^n$ and $\bar A = \la A_\sigma ;\; \sigma
\in\omega^{<\omega} \ra \leq \bar B$ such that for 
every $\sigma \in\omega^{<\omega}$
\sm
\ce{$T_\sigma(\bar A) \forces_\LA " \forall n \; (\breve f
\restrict I_n \in \phi_\sigma (n))".$ }}
\sm
{\capit Proof.} We construct $A_\sigma$ and $\phi_\sigma$
by recursion on $lh(\sigma)$. \par
$\tr$ $\sigma =\la\ra$. Using a standard fusion argument
for Laver forcing, find $T \leq_0 T_{\la\ra} (\bar B)$ and $\phi_{
\la\ra}$ so that
\sm
\ce{$T \forces_\LA " \forall n \; (\breve f \restrict I_n 
\in \phi_{\la\ra} (n) ) "$.}
\sm
\no For $\tau\in T$, let $A_\tau = succ_T (\tau)$ and
$\phi_\tau = \phi_{\la\ra}$.
\par
$\tr$ $lh(\sigma) > 0$. Either $A_\sigma$ and $\phi_\sigma$
have been constructed already, or we proceed as before.
$\qed$
\bigskip
Next, following Goldstern, Repick\'y, Shelah and Spinas
[GRSS, section 1], for $\bar A, \bar B \in Q$ we write
$\bar A \leq^* \bar B$ iff $\forall \sigma \; (A_\sigma
\sub^* B_\sigma )$ and $\forall^\infty \sigma \;
(A_\sigma \sub B_\sigma )$. Then:
\sm
{\bolds 3.6. Lemma.} (Judah -- Miller -- Shelah)
[JMS, section 3] (see also
[GRSS, section 1]) {\it Assume $MA_\lambda
(\sigma$--centered$)$. If $\la \bar A^\alpha ;\; \alpha < 
\lambda \ra$ is a $\leq^*$--descending sequence of elements
of $Q$, then there exists $\bar A^\alpha$ such that
for all $\alpha < \lambda$, $\bar A^\lambda \leq^* \bar A^\alpha$.
$\qed$}
\bigskip
We now have what we need to show that both $unif({\cal L})$
and $unif({\cal M})$ equal $2^\omega = \kappa$ after one Laver
real over a model of $MA(\sigma$--centered). Assume $\lambda
< \kappa$ and $\la \breve f_\alpha ;\; \alpha <\lambda\ra$
is a sequence of $\LA$--names for objects in $2^\omega$.
Given $T\in\LA$ we shall produce $S\leq T$ and
a set $D \sub 2^\omega$ which is a union of countably many
closed measure zero sets so that
\sm
\ce{$S\forces_\LA "\forall \alpha <\lambda \; (\breve f_\alpha
\in D)"$.}
\sm
\no 
Without loss $T=\omega^{<\omega}$. Using Lemmata 3.5. and 3.6.
construct a $\leq^*$--decreasing sequence $\la \bar A^\alpha
;\; \alpha \leq \lambda\ra$ of elements of $Q$ and $\la
\phi^\alpha_\sigma\in\prod_n [2^{I_n}]^n ; \;\alpha < \lambda ,
\sigma \in\omega^{<\omega}\ra$ so that for all $\alpha < \lambda$
and $\sigma \in\omega^{<\omega}$,
\sm
\ce{$T_\sigma(\bar A^\alpha) \forces_\LA "\forall n \;
( \breve f_\alpha \restrict I_n \in \phi^\alpha_\sigma (n))"$.}
\sm
\no Set $D^\alpha_\sigma := \{ f\in 2^\omega ;\; \forall n
\; (f\restrict I_n \in \phi^\alpha_\sigma (n)) \}$;
$D^\alpha_\sigma$ is a closed null set. Let $D: = \bigcup_{
\alpha < \lambda, \sigma \in\omega^{<\omega}} D^\alpha_\sigma$.
It is easy to check that 
\sm
\ce{$T_{\la\ra} (\bar A^\lambda ) \forces_\LA " \forall \alpha
< \lambda \; (\breve f_\alpha \in D)"$.}
\sm
\no By a recent result of Bartoszy\'nski and Shelah [BS],
the additivity of the $\sigma$--ideal generated by the
closed null sets equals $add({\cal M})$ and is thus $\geq 
{\bf p}$. Hence the assumptions on $V$ 
imply that $D$ is indeed contained in the
union of countably many closed null sets.
\sm
(c) It remains to show that $cov({\cal L}) = \omega_1$ in
the model of part (c) of Theorem 3.4. The proof is modelled
on the proof of [Br 1, 1.1., Lemma 1]. Let $G_\kappa$ be 
$\PP_\kappa$--generic over $V$, and work in $V[G_\kappa]$.
Suppose $\breve f$ is an $\LA$--name for an element of the
Cantor space $2^\omega$. Given $T\in\LA$ we shall find 
$S\leq T$ and a Borel null set $A\in V$ so that 
\sm
\ce{$(*)$ \hskip 2truecm $S\forces_\LA "\breve f \in A "$.}
\sm
\no Since $V \models GCH$, this is clearly sufficient.
\par
Without loss, $T = \omega^{<\omega}$. Fix $n\in\omega$.
Find $\bar A^n = \la A^n_\sigma ;\; \sigma\in\omega^{<\omega}\ra
\in Q$ and $t^n_\sigma \in 2^n $ ($\sigma \in \omega^{<\omega}$)
so that
\sm
\ce{$T_\sigma (\bar A^n) \forces_\LA " \breve f \restrict n=t^n_\sigma
"$}
\sm
\no for all $\sigma\in\omega^{<\omega}$.
$\bar A^n$ and $t^n_\sigma$ are easily constructed by recursion
on $lh(\sigma)$ (and we refer the reader to similar constructions
in 3.5. and 5.2. in case of doubt). Unfix $n$ and
find $\bar A \in Q$ so that $\bar A \leq^* \bar A^n$ for all
$n\in\omega$ (cf 3.6.). Using K\"onig's Lemma for $\omega$--trees
find $x^n_\sigma \geq n$ so that $m < n$ implies
$t^{x_\sigma^m}_\sigma \restrict m \sub t_\sigma^{x^n_\sigma}
\restrict n$. Define $f_\sigma \in 2^\omega$
by $f_\sigma (n) : = t^{x^{n+1}_\sigma}_\sigma (n)$.
\par
Because $\PP_\kappa$ is a finite--support iteration of
$\sigma$--centered p.o.'s, there are no reals random
over $V$ in $V[G_\kappa]$ (this was proved by Miller, and
is implicit in [Mi 1, section 5]; for an explicit argument,
see, e.g., [Br 1, subsection 1.1.]).
Thus there are measure--zero sets $A_\sigma \in V$ ($\sigma
\in\omega^{<\omega}$) so that $f_\sigma \in A_\sigma$.
By $ccc$--ness of $\PP_\kappa$, the union of the $A_\sigma$
(in $V[G_\kappa]$) is contained in a countable union of null
sets in $V$; call this union $A$; clearly $A$ is null. Working
still in $V$ find compact trees $T_n$ ($n\in\omega$) so that $A
\sub 2^\omega \setminus \bigcup_n [T_n]$.
We claim that (in $V[G_\kappa]$)
\sm
\ce{$T_{\la\ra} (\bar A) \forces_\LA " \breve f \in 2^\omega
\setminus \bigcup_n [T_n] "$,}
\sm
\no thus proving $(*)$.
\par
For suppose not. Find $S\leq T_{\la\ra} (\bar A)$ and $n\in\omega$
so that $S \forces_\LA " \breve f \in [T_n]"$. Let $\sigma :=
stem (S)$. As $f_\sigma \not\in [T_n]$, find $m\in\omega$ so that
$t_\sigma^{x_{\sigma}^m } \restrict m = f_\sigma \restrict m
\not\in T_n$; i.e. $t_\sigma^{x_\sigma^m} \not\in T_n$.
Clearly $S \leq T_\sigma (\bar A)$, and, throwing out finitely
many sequences (and their extensions), if necessary, we may
assume $S \leq T_\sigma (\bar A^{x^m_\sigma})$.
Thus
\sm
\ce{$S \forces_\LA " \breve f \restrict x^m_\sigma = t_\sigma^{
x^m_\sigma} \not\in T_n "$,}
\no a contradiction.
\par
This concludes the proof of Theorem 3.4. $\qed$
\bigskip
We continue with some comments on 3.4.: 
by (c) it is consistent that $cov ({\cal L}) = \omega_1$
after one Laver real over a model for $MA(\sigma$--centered).
However, we do not know whether this is always the case.
We think, though, that the following question has a
positive answer.
\sm
{\bolds 3.7. Question.} {\it Assume $V \models " MA + 2^\omega
=\kappa "$. Does $V^\LA \models " cov ({\cal L}) = \kappa "$?}
\sm
\no We shall comment on this again at the end of section 5.
Concerning the
cardinals which do not figure in Cicho\'n's diagram we get
${\bf p} = {\bf h} = {\bf s} = \omega_1$ 
after adding a Laver real. 
Closely related to Question 3.7. is:
\sm
{\bolds 3.8. Question.} {\it Assume $V \models " MA
(\sigma$--centered$) + 2^\omega = \kappa "$. Does 
$V^\LA \models " {\bf r} = 2^\omega "$?}
\sm
\no Finally
we do not
know the answer to the following.
\sm
{\bolds 3.9. Question.} {\it Is ${\bf a} = \omega_1$ after one
Laver real?}
\sm
\no We conjecture the answer is yes; in fact, it is an open
problem whether ${\bf d} < {\bf a}$ is consistent at all.
\par
Classically 
the consistency
of $unif({\cal L}) = unif({\cal M}) = 2^\omega > {\bf d}
= cov ({\cal L}) = \omega_1$ was shown by Bartoszy\'nski,
Judah and Shelah [BaJS, section 3] by iterating Shelah
($\QQ_{f,g}$) reals and infinitely often equal reals.
This proof works only in case $2^\omega = \omega_2$. In so 
far our result in Theorem 3.4. (c) is more general --- it
allows the continuum to be of arbitrary regular uncountable 
size. This partially answers an instance of Question 1 
in [Br 1]. We believe it is an interesting line of research to
figure
out whether other consistency results in Cicho\'n's
diagram can be gotten in a
similar fashion.
\bigskip
We next discuss several other forcing notions introduced
in section 1.
\sm
{\bolds 3.10. Corollary.} {\it (a) Let $\LA_{\cal U}$ be
Laver forcing with an ultrafilter ${\cal U}$. Then $V^{\LA_{\cal U}}
\models {\bf d} = \omega_1$.
\par
(b) Let $\MM_{\cal U}$ be Mathias forcing with a Ramsey
ultrafilter ${\cal U}$. Then $V^{\MM_{\cal U}} \models {\bf d}= \omega_1$.
\par
(c) Let $\MM$ be Mathias forcing. Then $V^\MM \models  {\bf d}
=\omega_1$.
\par
(d) Both $\LA \star \breve \LA$ and $\MM \star \breve \MM$ (the
two--step iterations) force $CH$ (and even $\diamondsuit$).}
\sm
{\capit Proof.} (a) $\LA_{\cal U}$ is trivially Laver--like and
closed under finite changes. Strong fusion is proved in a similar
fashion as for $\LA$ (see the proof of Theorem 3.4., part (a);
or see the original proof in [JS 1, Lemma 1.6.]).
\par
(b) Immediate from $\MM_{\cal U} \cong \LA_{\cal U}$ (see 1.1.).
\par
(c) Immediate from $\MM \cong P(\omega) / fin \star \MM_{\breve 
{\cal U}}$
(see 1.1.).
\par
(d) $CH$ is immediate from 3.4. (a) ((c), respectively) and the fact that
$\LA$ collapses the continuum below ${\bf h} \leq {\bf d}$
[GRSS, Theorem 2.7.] (that Mathias forcing collapses the
continuum onto ${\bf h}$, respectively).
\par
In the absence of $CH$ (in the ground model), $\diamondsuit$
follows from a result of Shelah's [RS]. If $CH$ holds,
however, then $\LA$ and $\MM$ can be shown to force
$\diamondsuit$ by an argument similar to the one for
Sacks forcing, due to Carlson and Laver [CL]. $\qed$
\sm
For Mathias forcing we can show similar results as in Theorem 3.4.;
more explicitly, after one Mathias real over a model for $MA
(\sigma$--centered),
$unif({\cal L}) = unif({\cal M}) = 2^\omega$. 
Questions 3.7., 3.8. and 3.9. are open for Mathias forcing
as well.
\bigskip
All forcing notions generated by a name for a dominating real
we have seen in sections 2 and 3 (as well as in [BJS])
add an unbounded family of reals of size $\omega_1$. Accordingly
we ask:
\sm
{\bolds 3.11. Question.} {\it Let $\PP$ be a proper forcing notion
generated
by a name for a dominating real. Does $V^\PP \models {\bf b}
=\omega_1$?}
\sm
\no An earlier version of this paper contained the same
question without the adjective "proper". Since then 
Saharon Shelah has noted that a slight modification
of the proofs in [SW] gives (non--proper) counterexamples.
\par Finally note that although both Theorem 2.1. and Theorem 3.1.
say
that adding a real can have strong combinatorial consequences and
collapse many invariants, they are mutually exclusive (as can be
seen from the effect on Cicho\'n's diagram): a Luzin set of size
$2^\omega$ entails that every dominating family has size $2^\omega$
as well.
\Bigskip
{\bold 4. Entr'acte --- Mathias--absoluteness and the
projective Ramsey--property}
\Smallskip
We leave our combinatorial considerations for a while
to deal with consequences of the main result of the
last section (for Mathias forcing) to descriptive set
theory.
\par
A subset $A\sub [\omega]^\omega$ is called {\it Ramsey}
iff there is $a \in [\omega]^\omega$ so that $[a]^\omega
\sub A$ or $[a]^\omega \sub [\omega]^\omega \setminus A$.
A classical result of Silver's says that analytic ($\Sigma^1_1$)
sets are Ramsey [Si]. Following Judah [Ju 2, section 2], given
a universe of set theory $V$ and a forcing notion $\PP\in V$,
we say that $V$ is {\it $\Sigma^1_n - \PP$--absolute}
iff for every $\Sigma^1_n$--sentence $\phi$ with parameters
in $V$, we have $V \models \phi$ iff $V^\PP \models \phi$.
By Shoenfield's absoluteness Lemma [Je 1, Theorem 98], $V$ is
always $\Sigma^1_2 - \PP$--absolute. Recently Halbeisen
[Ha] proved that $\Sigma^1_3  - \MM$--absoluteness is
equivalent to {\it all $\Sigma^1_2$--sets are Ramsey}.
We shall show that at least one implication is still true
one level higher in the projective hierarchy.
\smallskip
{\slan 4.1. Theorem.} {\it $\Sigma^1_4 - \MM$--absoluteness
implies that all $\Sigma^1_3$--sets are Ramsey.}
\smallskip
Roughly our argument follows the lines of the proof that
$\Sigma^1_4$--Amoeba--absoluteness implies
$\Sigma^1_3$--measurability [Br 2, section 2] --- still we are
indebted to Haim Judah who taught us many of the arguments 
involved and to Lorenz Halbeisen from whose talk at the Berlin
conference on Cantor's Set Theory we learned as well.
Before starting the proof itself, we shall give some
auxiliary results which might be of own interest.
\bigskip
{\slan 4.2. Proposition.} {\it Assume $\omega_1^{L[r]}
=\omega_1^V$, where $r\in\omega^\omega$. Then the reals of
$L[r][m]$ dominate the reals of $V[m]$, where $m$ is
Mathias ($\MM$--generic) over $V$.}
\smallskip
{\capit Remark.} Note that, in general, $m$ is {\it not}
Mathias over $L[r]$.
\smallskip
{\capit Proof.} Just note that the $\MM$--names $\breve f_\alpha$
($\alpha < \omega_1$) set up at the beginning of the proof
of Theorem 3.1. can all be defined within $L[r]$.
Therefore their interpretation $\breve f_\alpha [m]$ lies
in $L[r][m]$. The proposition now follows from the fact
that they form a dominating family in $V[m]$. $\qed$
\bigskip
{\slan 4.3. Proposition.} {\it $\Sigma^1_4 -\MM$--absoluteness
implies that $\omega_1^{L[r]} < \omega_1^V$ for all reals $r$.}
\sm
{\capit Proof.} The proof follows the lines of the proof of
[BJS, Theorem 2.6.]. Suppose, by way of contradiction, that
there is an $r\in\omega^\omega$ so that $\omega_1^{L[r]}
=\omega_1^V$. $\Sigma^1_3 -\MM$--absoluteness implies
--- by Halbeisen's result quoted at the beginning of
this section --- that the $\Sigma^1_2$--Ramsey property holds:
thus, by [Ju 1, Theorem 0.9.(b)], there are dominating
reals over all $L[s]$, where $s\in\omega^\omega$.
I.e.
\sm
\ce{$V\models \forall s \in\omega^\omega \; \exists
t\in\omega^\omega \; \forall u \in\omega^\omega \; ( u
\not\in L[s] \;\lor\; \exists n\;\forall m \geq n \;
(t(m) \geq u(m))),$}
\sm
\no which is a $\Pi^1_4$--statement. By $\Sigma^1_4
-\MM$--absoluteness, it should hold in $V[m]$. By the 
previous Proposition, however, it is false for $s =
\la r,m \ra$, a contradiction. $\qed$
\bigskip
In a sense the proof of 4.1. splits into two parts.
the first is 4.3., the second, given below, is the argument
showing that $\Sigma^1_4 - \MM$--absoluteness $+ \;\forall
r \; (\omega_1^{L[r]} < \omega_1 )$ entails the
$\Sigma^1_3$--Ramsey property. After having figured out
the latter part (as presented below), Haim Judah informed
us that he had proved that independently (as a consequence
of a general result about Souslin proper forcing notions).
\sm
{\capit Proof of Theorem 4.1.} Let $\phi (x)$ be a 
$\Sigma^1_3$--formula. First assume that $V[m]
\models \neg \phi (m)$, where $m$ is $\MM$--generic
over $V$. By a well--known result of Mathias [Ma, 2.9.
and 4.12.],
there is $S \in [\omega]^\omega$ so that 
$(\la\ra , S) \forces_\MM \neg \phi (\breve m)$, and
$m \sub S$. As subsets of Mathias--generics are
again Mathias--generic [Ma, 4.10.], and $m' \sub m$
implies (trivially) $m' \sub S$, we must have that
$V[m] \models \forall m' \in [m]^\omega \;
(\neg \phi (m'))$; i.e. $V[m] \models \exists x \;
\forall y \in [x]^\omega \; (\neg \phi (y))$.
By $\Sigma^1_4 - \MM$--absoluteness, this must hold in $V$ as
well, and we are done.
\par
So assume $V[m] \models \phi (m)$. Write $\MM =
P(\omega) / fin \star \MM_{\breve {\cal U}}$, where
$\breve{\cal U}$ is the $P(\omega) / fin$--name
for the generic Ramsey--ultrafilter (see 1.1.). Now note
that given $S\in P(\omega) / fin$ and $r\in\omega^\omega$
there is $T \sub^* S$ so that $T \forces_{P(\omega)
/fin} " \breve {\cal U} \cap L[r][S]$ is $P(\omega)
/fin$--generic over $L[r][S]$". To see this recall
that by the previous Proposition, $(P(\omega) / fin)^{L[r][S]}$
is countable (in the sense of $V$), and only countably many
maximal antichains of $(P(\omega) / fin)^{L[r][S]}$
lie in $L[r][S]$. Hence we can find $T \sub^* S$ which is
almost included in one member of each of the antichains.
In fact, we can now define $\bar{\cal U} :=
\{ S' \in P(\omega)  \cap L[r][S] ;\; T \sub^* S' \}$
in $V$, and $\bar{\cal U}$ is a Ramsey ultrafilter in
$L[r][S][\bar{\cal U}]$ (by genericity there are no new
reals in $L[r][S][\bar{\cal U}]$).
\par
Let $\phi (x) = \exists y\; \psi (y,x)$, where $\psi$ is
$\Pi^1_2$. Then $V[m] \models \psi (\breve s [m] , m)$
for some $\MM_{{\cal U}}$--name $\breve s$, where ${\cal U}$
is the generic Ramsey ultrafilter. Let $r$ code $\breve s$
as well as the parameters of $\psi$. By the previous 
paragraph and a density argument we can find $S \in
P(\omega) / fin \cap V$ so that $\bar{\cal U} :=
{\cal U} \cap L[r][S]$ is a Ramsey ultrafilter in
$L[r][S][\bar{\cal U}]$ (and $S \in {\cal U}$).
By a characterization of $\MM_{\cal U}$--genericity due
to Mathias ([Ma, 2.0.], see also [JS 1, Theorem 1.19]),
$m$ must be $\MM_{\bar{\cal U}}$--generic over $L[r]
[S][\bar{\cal U}]$. By Shoenfield's absoluteness Lemma,
$L[r][S][\bar{\cal U}][m] \models \psi (\breve s[m] , m)$.
Hence there is $T \in \bar{\cal U}$ so that
\sm
\ce{$L[r][S][\bar{\cal U}] \models \; (\la\ra , T)
\forces_{\MM_{\bar{\cal U}}} \psi (\breve s, \breve m).$}
\sm
\no Recall that $L[r][S][\bar{\cal U}] \subseteq V$, and
$\bar{\cal U}$ is a real in the sense of $V$; so the
previous Proposition applies again, and we get an
$\MM_{\bar{\cal U}}$--generic (over $L[r][S][\bar{\cal U}]$) 
$m$ in $V$ with $m \sub T$; thus
\sm
\ce{$V \models \psi(\breve s [m] , m)$}
\sm
\no and for all $m' \sub m$, $V \models \psi (\breve s [m'] ,
m')$ (by the argument of the first paragraph:
subsets of generics are generic). Hence
$V \models \forall m' \in [m]^\omega \; (\phi (m'))$, and
we are done. $\qed$

\Bigskip
%

{\bold 5. Counterpoint --- preserving MA($\sigma$--centered)}
\Smallskip
The following result shows that Miller's rational perfect set
forcing $\MI$ behaves rather differently from Laver and
Mathias forcings.
\sm
{\slan 5.1. Theorem.} {\it Miller forcing preserves the
axiom $MA(\sigma$--centered$)$; i.e. if $\kappa$ is regular
uncountable and $V \models " 2^\omega = \kappa \;\land\;
MA(\sigma$--centered$)"$, then $V^\MI \models "
2^\omega = \kappa \;\land\; MA(\sigma$--centered$)"$.}
\sm
{\capit Proof.} It was proved by Judah, Miller and Shelah in
[JMS, Theorem 4.1.] that $\MI$ doesn't collapse cardinals under
the assumptions. One main ingredient of our argument is the method
of their proof, as corrected by Goldstern, Johnson and Spinas
in [GJS, section 4] (note the similarity between
some arguments here and in the proof
of Theorem 3.4.). We start with setting up some notation
which is close to the one in the latter work. The proof then
will be broken up into a series of lemmata.
\par
Following [GJS, 4.3.] we call a sequence $\bar P = \la
P_\sigma ;\; \sigma\in\omega^{<\omega} \ra$ {\it good}
iff each $P_\sigma \sub \omega^{<\omega}$ is infinite,
$\tau \in P_\sigma$ implies $\sigma \subset \tau$, and
for $\sigma \in \omega^n$, if $\tau , \tau ' \in P_\sigma$
and $\tau \neq \tau '$, then $\tau (n) \neq \tau ' (n)$.
With $\sigma \in \omega^{<\omega}$ and a good $\bar P
=\la P_\sigma ;\; \sigma \in \omega^{<\omega} \ra$ we can
associate a condition $T_\sigma = T_\sigma (\bar P)  \in
\MI$ as follows: let $S$ be the smallest subset of
$\omega^{<\omega}$ such that $\sigma \in S$ and if $\tau \in
S$ then $P_\tau \sub S$; $T_\sigma$ is the unique condition
with $S = split (T_\sigma)$; i.e., $\sigma = stem (T_\sigma)$
and if $\tau \in split (T_\sigma)$, then $succ_{T_\sigma} (\tau)
= P_\tau$. We put $\bar P = \la P_\sigma ;\; \sigma \in
\omega^{<\omega} \ra \leq \bar Q = \la Q_\sigma ;\; \sigma \in
\omega^{<\omega} \ra$ iff $T_\sigma (\bar P) \leq T_\sigma
(\bar Q)$ for all $\sigma \in \omega^{<\omega}$ iff
$P_\sigma \sub split (T_\sigma (\bar Q))$ for all $\sigma
\in \omega^{<\omega}$ (see [GJS, 4.5.]; we note that our
relation $"\leq "$ corresponds to their $"\geq "$, because
our forcing--theoretic notation goes the other way round).
Finally, given a $\MI$--name $\breve A$ for an infinite
subset of $\omega$, we say a good $\bar P =\la P_\sigma
;\; \sigma \in\omega^{<\omega} \ra$ is {\it $\breve A$--nice}
iff
\sm
\item{$(\spadesuit)$} whenever $\sigma \in \omega^{<\omega}$ and
$\tau \in P_\sigma$, then $T_\tau (\bar P)$ decides the 
first $lh (\sigma) + \tau (lh (\sigma))$ values of
$\breve A$ (more explicitly, there is $m > lh (\sigma)
$ and $a_\tau = a \in [m]^{lh(\sigma) + \tau(lh(\sigma))}$
so that $T_\tau (\bar P) \forces_\MI " \breve A \cap m = a "$).
\bigskip
{\slan 5.2. Lemma.} {\it Given a $\MI$--name $\breve A$ for an
element of $[\omega]^\omega$, and a good $\bar P = \la
P_\sigma ;\; \sigma \in \omega^{<\omega} \ra$, there is
$\bar Q = \la Q_\sigma ;\; \sigma \in \omega^{<\omega} \ra
\leq \bar P$ which is good and $\breve A$--nice.}
\sm
{\capit Proof.} By recursion on $lh(\sigma)$ we construct
$Q_\sigma$ and an auxiliary $\tilde Q_\sigma$.
\par
$\tr$ $\sigma = \la\ra$. We let $\tilde Q_{\la\ra}
=P_{\la\ra}$, and for each $\tau \in P_{\la\ra}$ we find
$T ' (\tau) \leq T_\tau (\bar P)$ deciding the first 
$\tau (0)$ values of $\breve A$. Let $Q_{\la\ra} =
\{ stem (T' (\tau)) ; \tau \in P_{\la\ra} \}$. Clearly
$Q_{\la\ra} \sub split (T_{\la\ra} (\bar P))$.
\par
$\tr$ $lh(\sigma) \geq 1$. If $\sigma \not\in Q_{\sigma \restrict
n}$ for any $n < lh(\sigma)$ we proceed as in the previous case
(and guarantee that $Q_\sigma \cap Q_{\sigma \restrict n}
=\emptyset$ for $n < lh(\sigma)$).
Otherwise let $m < lh (\sigma)$ be unique with $\sigma
\in Q_{\sigma \restrict m}$. In this case $\sigma = stem
(T'(\tilde\tau)) $ (where the tree $T' (\tilde\tau)$
was constructed as a subtree of a tree with stem $\tilde\tau
\in \tilde Q_{\sigma\restrict m}$ when defining $Q_{\sigma \restrict
m}$). We let $\tilde Q_\sigma = succ_{T'(\tilde\tau)} (\sigma)$,
and for each $\tau \in \tilde Q_\sigma$ we find $T' (\tau)
\leq (T' (\tilde\tau))_\tau$ deciding the first $lh (\sigma)
+ \tau (lh (\sigma))$ values of $\breve A$. Let $Q_\sigma =
\{ stem (T' (\tau)) ; \; \tau \in \tilde Q_\sigma \}$.
By replacing $Q_\sigma$ by a cofinite subset, if
necessary, we can assume that $Q_\sigma \cap Q_{\sigma \restrict
n} = \emptyset$ for $n < lh (\sigma)$.
Clearly $Q_\sigma \sub split (T_\sigma (\bar P))$.
\par
This concludes the definition of the $Q_\sigma$. $\bar Q \leq
\bar P$ as well as the goodness of $\bar Q$ are immediate. It
follows
easily from the construction that $T_\tau (\bar Q) \leq
T' (\tilde\tau)$ for $\tau \in Q_\sigma$, where $\tilde\tau
\in \tilde Q_\sigma$ is such that $\sigma \subset
\tilde\tau \sub \tau$, and thus the former condition indeed satisfies
$(\spadesuit)$. $\qed$
\bigskip
Next, following again [GJS, section 4], for $\bar P = \la
P_\sigma ;\; \sigma \in \omega^{<\omega} \ra$ and $\bar Q
=\la Q_\sigma ; \; \sigma \in \omega^{<\omega} \ra$ good, we
say $\bar P \approx \bar Q$ iff $\forall \sigma \;
(P_\sigma =^* Q_\sigma )$ and $\forall^\infty \sigma \;
(P_\sigma = Q_\sigma )$; and we write $\bar P \leq^* \bar Q$
iff there is $\bar P ' \approx \bar P$ so that $\bar P '
\leq \bar Q$. Then:
\sm
{\slan 5.3. Lemma.} (Goldstern -- Johnson -- Spinas) [GJS,
section 4] \par {\it (a) $\leq^*$ is transitive.
\par
(b) Assume $MA_\lambda (\sigma$--centered$)$. If $\la \bar 
P^\alpha ;\; \alpha <\lambda\ra$ as a $\leq^*$--decreasing
sequence of good sequences, then there exists $\bar P^\lambda$
such that for all $\alpha < \lambda$, $\bar P^\lambda
\leq^* \bar P^\alpha$.} $\qed$
\bigskip
Now we have the main tools to be able to show the preservation
of $MA(\sigma$--centered$)$. By Bell's Theorem (cf 1.2.)
it suffices to show that $V^\MI \models "{\bf p} = \kappa"$.
To do this let $\lambda < \kappa$, and assume that $\la \breve
A_\alpha ;\; \alpha < \lambda \ra$ is a sequence of $\MI$--names
for infinite subsets of $\omega$ so that 
\sm
\ce{$\forces_\MI "\la \breve A_\alpha ;\; \alpha < \lambda \ra$
has the strong finite intersection property ".}
\sm
\no We have to prove that, given $T \in\MI$ there are $S\leq T$
and a $\MI$--name $\breve A$ for a subset of $\omega$ so that
\sm
\ce{$(+)$ \hskip 2truecm 
$S \forces_\MI " \vert \breve A \vert = \omega \;\land\;
\forall \alpha < \lambda \; (\breve A \sub^* \breve A_\alpha)"$.}
\sm
\no Without loss $T = \omega^{<\omega}$. Let $\la \Gamma_\alpha
; \; \alpha < \lambda \ra$ be an enumeration of the finite subsets
of $\lambda$. With each $\Gamma_\alpha$ we can associate a
$\MI$--name $\breve A_{\Gamma_\alpha}$ for an infinite
subset of $\omega$ so that 
\sm
\ce{$\forces_\MI " \breve A_{\Gamma_\alpha} = \bigcap_{\beta
\in \Gamma_\alpha} \breve A_\beta ".$}
\sm
\no Using Lemmata 5.2. and 5.3. (b) we can construct a
$\leq^*$--decreasing sequence $\la \bar P^\alpha ; \; \alpha
\leq \lambda \ra$ of good sequences so that $\bar P^\alpha$
is $\breve A_{\Gamma_\alpha}$--nice for $\alpha < \lambda$.
Let $S' := T_{\la\ra} (\bar P^\lambda)$. 
We shall define the p.o. $\PP$ for shooting a $\MI$--name
for a subset of $\omega$ through $S '$.
\par
For $\sigma \in split (S')$, we let $pred (\sigma)$ be the 
predecessor of $\sigma$ in $split (S')$. By $\breve
A_{\Gamma_\alpha}$--niceness of $\bar P^\alpha$ we have, for
$\sigma \in \omega^{<\omega}$ and $\tau \in P^\alpha_\sigma$,
$a^*_{\alpha ,\tau} \in 2^{\geq lh(\sigma) 
+ \tau (lh (\sigma))}$ so that
\sm
\ce{$T_\tau (\bar P^\alpha ) \forces_\MI " \breve A_{\Gamma_\alpha}^*
\restrict lh (a_{\alpha , \tau}^* ) = a^*_{\alpha, \tau} "$}
\sm
\no (where, for $a \in [\omega]^{\leq \omega}$, $a^* \in 2^{\leq
\omega}$ is the characteristic function of $a$ --- or of the
restriction of $a$ to $lh (a^*)$).
We say that $\sigma \in split(S')$ is {\it $\Gamma_\alpha$--happy} iff
$T_\sigma (\bar P^\lambda) \leq T_\sigma (\bar P^\alpha )$. We note
that for $\alpha < \lambda$, by definition of $\leq^*$ and the
construction, almost all $\sigma \in split (S')$ are
$\Gamma_\alpha$--happy.
We define:
\medskip
\ce{$(\Sigma , (a^*_\sigma ;\;\sigma\in\Sigma ) , \Gamma ) \in \PP
\Longleftrightarrow \Sigma \sub split (S') $ is finite and closed
under
predecessors in $split(S') \;\land$}
\par
\ce{$\land\; a^*_\sigma \in 2^{<\omega} \;\land\; lh (a^*_\sigma)
\geq lh (pred(\sigma)) + \sigma(lh(pred(\sigma))) \;\land\;
(\sigma \sub\tau \to a^*_\sigma \sub a^*_\tau ) \;\land\;
\Gamma \sub \lambda$ is finite}
\medskip\sm
\ce{$(\Sigma ' , (a^{*\prime}_\sigma ;\; \sigma\in\Sigma ') ,
\Gamma ' ) \leq (\Sigma , (a^*_\sigma ;\;\sigma\in\Sigma ) , 
\Gamma ) \Longleftrightarrow \Sigma' \supseteq \Sigma \;\land\;
(\sigma \in\Sigma \to a^{*\prime}_\sigma = a^*_\sigma ) \;\land\;
\Gamma ' \supseteq \Gamma \;\land$}
\par
\ce{$\land\; \forall\sigma\in\Sigma ' \setminus\Sigma
\;\forall\alpha
\in\Gamma\; ($ if $\sigma$ is $\{\alpha\}$--happy,}
\par
\ce{then $a_\sigma '
\cap [lh(a^{*\prime}_{pred(\sigma)}) , lh (a^{*\prime}_\sigma)) \sub
a_{\beta,\sigma} \cap [lh(a^{*\prime}_{pred(\sigma)}), lh (a^{*
\prime}_\sigma))$}
\par
\ce{
for some $\Gamma_\beta \supseteq \{\alpha\}$ so that
$\sigma$ is $\Gamma_\beta$--happy).}
\medskip
\no $\PP$ is easily seen to be $\sigma$--centered.
\sm
{\slan 5.4. Observation.} {\it Let $G$ be $\PP$--generic over
$N \prec H(\chi)$, where $\lambda \sub N$, $\PP \in N$, and $
\vert N \vert = \lambda$.
\par
(i) For all $\alpha < \lambda$ and almost all $\sigma
\in split (S')$, $a_\sigma \cap [lh (a^*_{pred (\sigma)})
, lh(a^*_\sigma) ) \sub a_{\beta,\tau} \cap [lh (a^*_{pred
(\sigma)}) , lh (a^*_\sigma))$ for any $\tau \supseteq
\sigma$ in $split (S')$ which is long enough (i.e. $lh (a^*_{
\beta, \tau}) \geq lh (a^*_\sigma)$), where $\Gamma_\beta
=\{ \alpha \}$.
\par
(ii) For all $\sigma \in split (S') \; \exists^\infty \tau
\in P^\lambda_\sigma \; (a_\tau \cap [lh (a^*_\sigma) , lh
(a^*_\tau)) \neq \emptyset )$.}
\sm
{\capit Proof.} (i) Use genericity and the fact that almost all
$\sigma \in split (S')$ are $\{\alpha\}$--happy.
\par
(ii) Fix $\sigma \in split (S')$; given $(\Sigma , (a^*_\tau
;\; \tau\in\Sigma ) , \Gamma) \in \PP$
with $\sigma\in\Sigma$, almost all $\tau \in 
P^\lambda_\sigma $ are 
$\Gamma$--happy and not yet in $\Sigma$. For such $\tau$
(with long enough $\tau(lh(\sigma))$) we can find an extension with
non--trivial intersection by construction. Thus genericity gives 
the desired result. $\qed$
\sm
Next define $S\leq S'$ by recursion on its levels. Put $\la\ra$
into $S$. Assume $\sigma \in S\cap split (S')$. Then put $\tau
\in P^\lambda_\sigma$
(as well as all $\tau \restrict n$) into $S$
iff $a_\tau \cap [lh (a^*_\sigma) , lh (a^*_\tau )) \neq \emptyset$.
By 5.4. (ii) $S$ is indeed a Miller tree.
Let $\breve A$ be the $\MI$--name defined by $S \forces_\MI
" \breve A = \bigcup_{n\in\omega\;{\rm so \; that}
\; \breve m \restrict n \in
split (S)} a_{\breve m \restrict n}"$, where $\breve m$ is the name
for the generic real.
It is easily seen that $S \forces_\MI " \vert \breve A \vert =
\omega "$.
\sm
{\slan 5.5. Claim.} {\it $S \forces_\MI "\forall \alpha < \lambda
\; (\breve A \sub^* \breve A_\alpha )"$.}
\sm
{\capit Proof.} Given $T \leq S$ and $\alpha < \lambda$ we have to
find $T' \leq T$ and $n\in\omega$ so that 
\sm
\ce{$(++)$ \hskip 2truecm $T' \forces_\MI " \breve A \setminus
n \sub \breve A_\alpha "$.}
\sm
\no To do this simply let $\sigma \in split (T)$ be so long that
for all $\tau \in split (T)$ extending $\sigma$ (i) of 5.4. is
satisfied. Next let $n = lh (a^*_{pred(\sigma)})$ and $T'
= T_\sigma$. $T'$ and $n$ are easily seen to satisfy $(++)$.
$\qed$
\sm
Thus $S$ satisfies $(+)$, and we are done with Theorem 5.1.
$\qed$
\bigskip
When studying the effect of a Miller real on the 
cardinals in Cicho\'n's diagram we have to assume
again $MA(\sigma$--centered) (or maybe something
slightly weaker) to avoid the pathology of collapsing cardinals.
Then we get:
\sm {\slan 5.6. Theorem.} {\it Let $V \models ZFC$. Then
\par
(a) $V^\MI \models " add({\cal L}) = \omega_1 "$.
\par
(b) If $V \models " MA(\sigma$--centered$) + 2^\omega = \kappa "$,
then $V^\MI \models " add ({\cal M}) = \kappa = 2^\omega "$.
\par
(c) If we force $MA(\sigma$--centered$)$ (as in Theorem 3.4. (c)),
and then add a Miller real, we have $add({\cal M}) = 2^\omega$
and $cov({\cal L})= \omega_1$.
\par
(d) If $V \models " MA + 2^\omega = \kappa "$, then $V^\MI
\models " cov ({\cal L}) = 2^\omega = \kappa "$.}
\sm
{\capit Proof.} 
(a) To see that $add({\cal L}) = \omega_1$ in $V^\MI$ (where $V$
is arbitrary), we use the combinatorial characterization of this
cardinal (see 1.2.). Let $\{ f_\alpha ;\; \alpha < \omega_1 \}
\in V$ be a family of almost--disjoint functions from $\omega$
to $\omega$ (cf the proof of Theorem 2.1.!). $\breve f_\alpha=
f_\alpha \circ \breve m$, where $\breve m$ is the $\MI$--name for
the generic, is a $\MI$--name for a new real. Let $\breve \phi$
be the $\MI$--name for a slalom (i.e. $\forces_\MI
"\breve \phi \in ([\omega]^{<\omega})^\omega \;\land\;
\forall n \; (\vert\breve \phi (n) \vert = n )"$) and assume that
\sm
\ce{$\forces_\MI " \forall \alpha < \omega_1 \; \forall^\infty n
\; (\breve f_\alpha (n) \in \breve\phi (n))"$.}
\sm
\no Let us say $\sigma \in T$ is an {\it $n$--th splitting node}
if $\sigma \in split (T)$ and for exactly $n$ many predecessors
$\tau\sub\sigma$, $\tau \in split (T)$. Now let $T$ be so
that: if $\sigma$ is an $n$--th splitting node, then $T_{\sigma
\hat{\;}\la m\ra} $ decides $\breve\phi (n)$ whenever
$\sigma\hat{\;}\la m \ra \in T$. Fix $\sigma \in split (T)$,
$lh(\sigma) = n$. Say $\alpha$ is {\it $\sigma$--bad} iff:
for almost all $\ell$ with $\sigma\hat{\;}\la\ell\ra\in T$
for all $\tau\in T$ with $\tau \supseteq \sigma\hat{\;}\la
\ell\ra$ and which is an $n$--th splitting node and for
(almost) all $m$ with $\tau\hat{\;}\la m \ra \in T$, the value
forced to $\breve\phi (n)$ contains $f_\alpha (\ell)$
[which is the value forced to $\breve f_\alpha (n)$ by
$T_{\sigma \hat{\;} \la \ell\ra}$].
At most $n$ many $\alpha$ can be $\sigma$--bad. Let $\alpha$ be
such that it is not $\sigma$--bad for any $\sigma \in split (T)$.
Construct $T' \leq T$ as follows. $stem (T') = stem (T) =:\sigma$;
$lh(\sigma) = n_{\la\ra}$. Choose infinitely many $\ell\in\omega$
so that $\sigma\hat{\;}\la \ell\ra \in T$ and there are $\tau_\ell
\in T$ with $\tau_\ell \supseteq \sigma\hat{\;}\la\ell\ra$ and
which are $n_{\la\ra}$--th splitting nodes so that for some
$m_\ell \in\omega$ with $\tau_\ell \hat{\;} \la m_\ell \ra
\in T$, the value forced to $\breve\phi (n_{\la\ra})$ doesn't
contain $f_\alpha (\ell)$. Let $\tilde \tau_\ell \in
split (T)$ extend $\tau_\ell \hat{\;} \la m_\ell\ra$.
We put all $\tilde\tau_\ell$ into $T'$. Repeat this procedure with
each $\tilde\tau_\ell$ in place of $\sigma$; etc. In the end we
construct a condition $T' \leq T$ so that 
\sm
\ce{$T' \forces_\MI "\exists^\infty n \; (\breve f_\alpha (n)
\not\in \breve\phi (n))"$.}
\sm
\no This concludes the proof of part (a).
\par
(b) Immediate from Theorem 5.1.
\par
(c) Rewrite the proof of Theorem 3.4. (c).
\par
(d) Let $\la s_n ;\; n\in\omega\ra$ enumerate the basic clopen sets
of $2^\omega$. Let $\breve O$ be a $\MI$--name for an open set.
Clearly we can identify $\breve O$ with a $\MI$--name $\breve f$ for
a strictly increasing function so that
\sm \ce{$\forces_\MI " \breve O = \bigcup_{n\in\omega}
s_{\breve f (n)} "$.}
\sm
\no We say a good $\bar P = \la P_\sigma ;\;
\sigma\in\omega^{<\omega}\ra$ is {\it $\breve O$--soft}
iff
\sm
\item{$(\heartsuit)$} there is a function $f : \omega^{<\omega}\to
\omega^{<\omega}$ so that $\forall \sigma, \tau \in 
\omega^{<\omega}$ \par
\itemitem{(I)} $\tau \in P_\sigma$ implies $f(\sigma) \sub
f (\tau)$; \par
\itemitem{(II)} $T_\tau (\bar P) \forces_\MI " f(\tau) = \breve f
\restrict lh (f (\tau)) "$; \par
\itemitem{(III)} $\forall k\in\omega\;\forall h\in [T_\tau 
(\bar P )] \;\exists \rho\in split (T_\tau (\bar P )) \;
(lh (f (\rho )) \geq k \;\land\; \rho \sub h )$; \par
\itemitem{(IV)} $\tau \in P_\sigma $ and $\tau (lh (\sigma ))
\geq m$ imply $T_\tau (\bar P) \forces_\MI " \mu
(\breve O \setminus \bigcup_{n \in lh (f(\tau ))}
s_{\breve f (n)} ) \leq {1 \over m} "$. \par
\bigskip
{\bolds 5.7. Lemma.} {\it Given $\MI$--names $\breve O$ and 
$\breve f$ as above and a good $\bar P = \la P_\sigma ;\; \sigma
\in\omega^{<\omega}\ra$, there is $\bar Q = \la Q_\sigma 
;\;\sigma\in\omega^{<\omega}\ra \leq \bar P$
which is good and $\breve O$--soft.}
\sm
{\capit Proof.} The proof is similar to the one of 5.2. ---
we proceed by recursion on $lh (\sigma)$, and construct
$Q_\sigma$, $\tilde Q_\sigma$ and $f$.
\par
$\tr$ $\sigma = \la\ra$. We let $\tilde Q_{\la\ra}
=P_{\la\ra}$. Let $f(\la\ra ) = \la\ra$. Fix $\tau \in
P_{\la\ra}$; find $T' (\tau ) \leq T_\tau (\bar P)$ and
$\ell \in\omega$, $\ell\geq 1$, so that 
\sm
\ce{$T' (\tau) \forces_\MI " \mu (\breve O \setminus \bigcup_{n
\in\ell} s_{\breve f (n) } ) \leq {1 \over \tau (0) }"$}
\sm
\no and $T' (\tau )$ decides the first $\ell$ values of $\breve f$;
say $T' (\tau ) \forces_\MI " \breve f \restrict \ell = \rho "$. 
Set $f (stem (T' (\tau ))) = \rho$, and --- unfixing $\tau
$ --- let $Q_{\la\ra} = \{ stem ((T ' (\tau )) ; \;
\tau \in P_{\la\ra} \}$.
\par
$\tr$ $lh (\sigma) \geq 1$. If $\sigma \not\in Q_{\sigma \restrict
n} $ for any $n < lh (\sigma)$ we proceed as above (letting
$f (\sigma ) = \la\ra$ and guaranteeing $Q_\sigma \cap Q_{
\sigma \restrict n} = \emptyset$ for $n < lh (\sigma )$).
Otherwise let $m < lh (\sigma )$ be unique with $\sigma \in
Q_{\sigma \restrict m}$. Then $\sigma = stem (T' (\tilde \tau ))$
(where $\sigma \restrict m \subset \tilde \tau \sub \sigma$
and $\tilde\tau \in \tilde Q_{\sigma \restrict m}$). Note
that $f(\sigma )$ has been defined already. Let $\tilde Q_\sigma
= succ_{T' (\tilde\tau )} (\sigma )$. Fix $\tau \in 
\tilde Q_\sigma$; find $T' (\tau ) \leq ( T' (\tilde\tau )
)_\tau$ and $\ell \geq max \{ lh (\sigma), lh (f (\sigma )) \}$
so that
\sm
\ce{$T' (\tau) \forces_\MI " \mu (\breve O \setminus 
\bigcup_{n\in\ell} s_{\breve f (n)} ) \leq { 1 \over \tau
(lh (\sigma )) } "$}
\sm
\no and $T' (\tau )$ decides the first $\ell$ values of $\breve f$;
say $T' (\tau ) \forces_\MI " \breve f \restrict \ell
= \rho "$. Set $f (stem (T' (\tau ))) = \rho$, and ---
unfixing $\tau $ --- let $Q_\sigma = \{ stem ( T' (\tau )) ;
\; \tau\in \tilde Q_\sigma \}$ (without loss disjoint
from $Q_{\sigma\restrict n}$ for $n\in lh (\sigma )$).
\par
Then $\bar Q \leq \bar P$ is easily seen to be good
and $\breve O$--soft. $\qed$
\sm
We note that if $\bar P \leq \bar Q$, $\bar P$ and $\bar Q$
are good, and $\bar Q$ is $\breve O$--soft, then $\bar P$ is
$\breve O$--soft as well, and this is witnessed
by the same function $f$. If $\bar P \leq^* \bar Q$, $\bar P$
and $\bar Q$ good and $\bar Q$ $\breve O$--soft, then
(I) -- (III) in $(\heartsuit)$ are satisfied (redefine
the function $f$ at finitely many places) and (IV) is
satisfied for almost all $\sigma$, and for $\sigma$ for
which it does not hold, it is still satisfied for almost
all $\tau \in P_\sigma$. We call $\bar P$ with these
properties {\it almost $\breve O$--soft}.
\bigskip
Fix $x\in 2^\omega$. We define the rank function $\rho_x =
\rho_x^{\bar P}$ on $\omega^{<\omega}$ (where ${\bar P}$
is good and almost $\breve O$--soft).
\sm
\ce{$\rho_x (\tau) = 0 \longleftrightarrow x \in 
\bigcup_{m < lh(f(\tau ))} s_{f(\tau) (m)}$;}
\sm
\ce{$\rho_x (\tau ) = \alpha \longleftrightarrow$ for no
$\beta < \alpha$ do we have $\rho_x (\tau) = \beta$}
\par
\ce{and $\forall^\infty \sigma \in P_\tau \; (\rho_x (\sigma )
< \alpha )$;}
\sm
\ce{$\rho_x (\tau) = \infty \longleftrightarrow$ for no
$\beta < \omega_1$ do we have $\rho_x (\tau) = \beta$.}
\sm
\no Note that the statement $" \rho_x (\sigma) = \infty "$
is equivalent to the existence of $\bar Q \leq \bar P$,
$\bar Q$ good, $\forall \tau \in \omega^{<\omega}
\; (Q_\tau \sub P_\tau)$ such that $\forall \tau
\in split (T_\sigma (\bar Q)) \; (x \not\in \bigcup_{
m < lh (f(\tau ))} s_{f(\tau )(m)})$. The latter is
easily seen to be a $\Sigma^1_1$--statement about
$x$ (with parameters $\bar P$ and $f$); hence $H_\sigma =
H_\sigma^{\bar P} := \{ x ;\;\rho_x (\sigma ) 
< \omega_1 \}$ is $\Pi^1_1$ and thus measurable.
\sm
{\bolds 5.8. Lemma.} {\it Let $\sigma \in\omega^{<\omega}$,
$\bar P$ good and almost $\breve O$--soft, and assume 
$T_\sigma (\bar P) \forces_\MI " \mu (\breve O) \leq \mu "$,
for some $\mu \leq 1$. Then $\mu (H_\sigma) \leq \mu$.}
\sm
{\capit Proof.} We show by induction on $\alpha < \omega_1$
simultaneously for all $\tau\in split (T_\sigma (\bar P))$
that $\mu (\{ x ;\; \rho_x (\tau ) \leq \alpha \}) \leq
\mu$.
\par
$\tr$ $\alpha = 0$. As $T_\tau (\bar P) \forces_\MI " \mu
(\breve O) \leq \mu "$ and $T_\tau (\bar P) \forces_\MI
"\bigcup_{m < lh(f(\tau ))} s_{f(\tau) (m)} \sub \breve O "$,
this is immediate.
\par
$\tr$ $\alpha > 0$. Note that $\{ x ;\; \rho_x (\tau ) \leq
\alpha \} = \bigcup_{i\in\omega} \bigcap_{\tilde\tau \in
P_\tau\;{\rm with}\;\tilde\tau (lh (\tau)) >i} \bigcup_{
\beta < \alpha} \{ x ; \; \rho_x (\tilde\tau )\leq\beta\}$.
The set on the right--hand side is easily seen to have
measure $\leq\mu$, by induction. $\qed$
\bigskip
Note that if $\bar Q \leq \bar P$, $\bar Q$ and $\bar P$ 
good and $\bar P$ almost $\breve O$--soft, and $\forall \sigma
\in\omega^{<\omega}\; (Q_\sigma \sub^* P_\sigma )$, then
$H_\sigma^{\bar Q} \supseteq H^{\bar P}_\sigma$ for 
all $\sigma\in\omega^{<\omega}$.
Let us say an almost $\breve O$--soft good $\bar P$ is 
{\it $\breve O$--maximal} iff whenever $\bar Q \leq^* \bar P$
good and $\forall\sigma\in\omega^{<\omega} \;
(Q_\sigma \sub^* P_\sigma )$, then $\forall \sigma\in\omega^{
<\omega} \; (\mu (H^{\bar Q}_\sigma ) = \mu (H^{\bar P}_\sigma
))$.
\sm
{\bolds 5.9. Observation.} {\it Given an almost $\breve O$--soft
good $\bar P$, there is $\bar Q \leq^* \bar P$ such that
$\forall \sigma\in\omega^{<\omega} \; (Q_\sigma \sub^* P_\sigma
)$ and $\bar Q$ is $\breve O$--maximal.}
\sm
{\capit Proof.} If $\bar P$ is almost $\breve O$--soft good and not
$\breve O$--maximal, then we can find $\bar Q \leq^* \bar P$
so that $\forall \sigma\in\omega^{<\omega} \; (Q_\sigma \sub^*
P_\sigma)$ and $\tau\in\omega^{<\omega}$ so that $\mu (H^{\bar
Q}_\tau) > \mu (H^{\bar P}_\tau )$. Thus if there were no
$\bar Q \leq^* \bar P$ with $\forall\sigma\in\omega^{<\omega}
\; (Q_\sigma \sub^* P_\sigma )$ which is $\breve O$--maximal,
we could construct a sequence $\la \bar P^\alpha ;\;\alpha
< \omega_1 \ra$ so that $\alpha < \beta$ implies $\bar P^\beta
\leq^* \bar P^\alpha$, $\forall\sigma\in\omega^{<\omega}
\; (P^\beta_\sigma \sub^* P^\alpha_\sigma )$ and $\forall
\alpha <\omega_1\;\exists\tau\in\omega^{<\omega}$
so that $\mu (H^{\bar P^{\alpha + 1}}_\tau) < \mu
(H^{\bar P^\alpha}_\tau )$, a contradiction. $\qed$
\bigskip
Equipped with these lemmata (as well as Lemma 5.3.)
we are ready for the proof that $cov({\cal L}) = \kappa$ in
the last model. Assume $\lambda < \kappa$ and $\la
\breve N^\alpha ;\;\alpha < \lambda\ra$ is a sequence
of $\MI$--names for null $G_\delta$'s. Let $\breve U^\alpha_n$,
$n \in\omega$, be $\MI$--names for open sets so that
\sm
\ce{$\forces_\MI " \forall n \; (\mu (\breve U^\alpha_n ) \leq
{1 \over 2^n }) \;\land\; \breve N^\alpha = \bigcap_n
\breve U^\alpha_n "$.}
\sm
\no Given $T\in\MI$ we shall produce $S \leq T$ and a $\MI$--name
$\breve x$ for an object in $2^\omega$ so that
\sm
\ce{$(*)$ \hskip 2truecm  $S \forces_\MI " \forall \alpha < \lambda
\; (\breve x \not\in \breve N^\alpha ) "$.}
\sm
\no Without loss $T = \omega^{<\omega}$. Let $\la \breve O^\alpha
;\;
\alpha < \lambda\ra$ be an enumeration of $\MI$--names for all
finite unions of the names $\breve U^\alpha_n$. Let
$\breve f^\alpha$, $\alpha < \lambda$, be the names for the
corresponding functions (i.e. $\forces_\MI " \breve O^\alpha
= \bigcup_n s_{\breve f^\alpha (n)} "$). Using Lemmata 5.7.
and 5.3. construct a $\leq^*$--decreasing sequence $\la 
\bar P^\alpha ;\;\alpha \leq \lambda\ra$ of good sequences so that for
all $\alpha < \lambda$, $\bar P^\alpha$ is $\breve O^\alpha$--soft.
Then $\bar P^\lambda$ is almost $\breve O^\alpha$--soft for all
$\alpha < \lambda$. Making another construction of length
$\lambda$, if necessary, using 5.9. and (essentially) 3.6.
(at limit steps), we can assume that $\bar P^\lambda$ is $\breve
O^\alpha$--maximal for all $\alpha < \lambda$. Let
$f^\alpha$ be the function witnessing this; and let $H^\alpha_\sigma
:= H^{\alpha , \bar P^\lambda}_\sigma$ be defined as above 
(before 5.8.) for $\alpha < \lambda$ and $\sigma\in\omega^{
<\omega}$.
\par
We define the p.o. $\QQ$ for shooting a Miller tree through
$S' := T_{\la\ra} (\bar P^\lambda)$ and producing a name
for a real. Conditions of $\QQ$ are of the form $p = (\Sigma^p,
(\bar B^p_\sigma ;\; \sigma \in \Sigma^p ), \Gamma^p )=
(\sigma , (\bar B_\sigma ;\;\sigma\in\Sigma ) , \Gamma )$ so that
\sm
\item{(i)} $\Sigma \sub split (S')$ is finite and closed 
under predecessors in $split (S')$; \par
\item{(ii)} $\bar B_\sigma = \la B_{\sigma , 0} , ... ,B_{
\sigma, n_\sigma - 1} \ra$ where $B_{\sigma , i} \in \BB$;
\par
\item{(iii)} $i<j<n_\sigma$ implies $B_{\sigma, i } \supseteq
B_{\sigma , j}$;
\par
\item{(iv)} $\tau\in P_\sigma \cap \Sigma$ and $\tau (lh (\sigma))
\geq i$ imply $i < n_\sigma - 1$ and $B_{\tau , j} \sub B_{\sigma ,
i}$ for any $j\in n_\tau$; \par
\item{(v)} if $\tau \in \Sigma$ is a final node (i.e. $P_\tau
\cap \Sigma = \emptyset$), then $n_\tau = 1$;
\par
\item{(vi)} $\Gamma \in [\lambda \times \omega]^{<\omega}$;
\par
\item{(vii)} letting $\alpha = \alpha (\Gamma)$ be such
that $\forces_\MI " \breve O^\alpha = \bigcup_{ (\beta,n)
\in \Gamma} \breve U^\beta_n "$, we have
$\forall \sigma\in\Sigma \; (H^\alpha_\sigma \cap B_{\sigma, n_\sigma
-1} = \emptyset)$.
\par\sm
\no We put $(\Sigma^p , (\bar B^p_\sigma ;\;\sigma\in\Sigma^p )
, \Gamma^p ) = p \leq q = (\Sigma^q , (\bar B^q_\sigma ; \; \sigma
\in\Sigma^q ) , \Gamma^q)$ iff
\sm
\item{(I)} $\Sigma^p \supseteq \Sigma^q$;
\par
\item{(II)} $n^p_\sigma \geq n^q_\sigma$ and $B^p_{\sigma ,i}
\subseteq B^q_{\sigma , i}$ for $\sigma \in \Sigma^q$
and $i < n^q_\sigma$; \par
\item{(III)} $\Gamma^p \supseteq \Gamma^q$.
\par\bigskip
{\bolds 5.10. Main Lemma.} {\it $\QQ$ is $ccc$.}
\sm
{\capit Proof.} Let $\la p^\beta = (\Sigma^\beta ,
(\bar B^\beta_\sigma ;\; \sigma\in\Sigma^\beta ) , \Gamma^\beta
) ;\; \beta < \omega_1 \ra$ be a sequence of elements of $\QQ$.
Going over to stronger conditions, if necessary, we can assume there
are $\Sigma$, $n_\sigma$ ($\sigma\in\Sigma$) and $B_{\sigma , i}$
($\sigma\in\Sigma , i\in n_\sigma$) so that  $\Sigma^\beta
=\Sigma$, $n_\sigma = n^\beta_\sigma$ and $B_{\sigma , i }
=B_{\sigma , i}^\beta$ for all $\beta < \omega_1$
(using $MA_{\omega_1}$ for the latter --- see [Tr 1,
Lemma 5.1.]). We have to find $\beta < \gamma$ so that,
letting $\alpha = \alpha (\Gamma^\beta \cup \Gamma^\gamma)$,
$\mu(B_{\sigma,n_\sigma - 1} \setminus H^\alpha_\sigma ) > 0$
for all $\sigma\in\Sigma$; then $p^\beta$ and $p^\gamma$
will be compatible. This is not trivial, for we may have
$\mu (H^\alpha_\sigma \setminus (H^{\alpha (\Gamma^\beta)}_\sigma
\cup H^{\alpha(\Gamma^\gamma)}_\sigma )) > 0$.
\par
Let $\Sigma = \{ \sigma_i ;\; i\in k\}$. We shall produce $\Delta_i
\in [\omega_1]^{\omega_1}$, $\Delta_{i+1} \sub \Delta_i$
so that when $\beta, \gamma \in \Delta_i$ then $\mu
(B_{\sigma_i , n_{\sigma_i} - 1} \setminus H^\alpha_{\sigma_i}
) > 0$ where $\alpha = \alpha (\Gamma^\beta \cup \Gamma^\gamma)$.
Fix $i<k$, and suppose $\Delta_{i-1}$ has been constructed 
(where $\Delta_{-1} = \omega_1$). Let 
$B:=B_{\sigma_i , n_{\sigma_i} - 1}$. Fix $\beta \in
\Delta_{i-1}$. We claim that there is $n^\beta$ so that
$\forall^\infty \tau \in P^\lambda_{\sigma_i} \;
(\mu (B \setminus H^{\alpha (\Gamma^\beta)}_\tau ) \geq
{1 \over n^\beta} )$. (**)
\par
For suppose not. Fix $g\in\omega^{<\omega}$ increasing and find
$\{ \tau_n ;\; n\in\omega \} \sub P^\lambda_{\sigma_i}$
so that $\mu (B \setminus H^{\alpha(\Gamma^\beta)}_{\tau_n}
) < {1 \over g(n)}$. Let $\bar P^* \leq \bar P^\lambda$
be defined by 
$$P^*_\sigma = \cases{ \{ \tau_n ; \; n\in\omega\},
& if $\sigma = \sigma_i$ \cr
P_\sigma, & otherwise. \cr}$$
Then $\mu (H^{\alpha (\Gamma^\beta), \bar P^*}_{\sigma_i}
\cap B) = \mu ((\bigcup_{m\in\omega} \bigcap_{n\geq m}
H^{\alpha(\Gamma^\beta)}_{\tau_n}
) \cap B) > 0$, contradicting $\breve O^{\alpha (\Gamma^\beta)}
$--maximality of $\bar P^\lambda$.
This shows (**).
\par
Without loss there is $n$ so that $n_\beta = n$ for all
$\beta \in \Delta_{i-1}$. Let $\Gamma \in [\Delta_{i-1}]^{
2n}$. Assume that for all pairs $\{ \beta , \gamma \} \in [
\Gamma ]^2$ we have $\mu (B \setminus H^{\alpha (\Gamma^\beta
\cup \Gamma^\gamma )}_{\sigma_i} ) = 0$ (***). Find
$\tau\in P^\lambda_{\sigma_i}$ so that $\tau (lh (\sigma_i))
\geq n^4$, $\tau$ satisfies (IV) in $(\heartsuit)$
for all $\beta \in\Gamma$, $\mu (B \setminus H^{\alpha (\Gamma^\beta
)}_\tau ) \geq {1 \over n}$ for all $\beta \in \Gamma$
(by (**)), and $\mu (B \setminus H_\tau^{\alpha (\Gamma^\beta
\cup \Gamma^\gamma)} ) < {1 \over n^4 }$ for all
$\{ \beta , \gamma \} \in [\Gamma]^2$ (by (***)).
Given $\beta \in \Gamma$, we have (by (IV) in $(\heartsuit)$)
\sm
\ce{$T_\tau (\bar P^\lambda) \forces_\MI " \mu (\breve O^{
\alpha (\Gamma^\beta)} \setminus \bigcup_{n \in lh (f^{
\alpha (\Gamma^\beta)} (\tau))} s_{\breve f^{\alpha (\Gamma^\beta
)} (n)} )
\leq {1 \over n^4} "$.}
\sm
\no Note that $\bigcup_{n \in lh (f^{\alpha (\Gamma^\beta)}
(\tau))} s_{f^{\alpha (\Gamma^\beta)} (\tau) (n)}
\sub H_\tau^{\alpha (\Gamma^\beta)}$; thus for $\{ \beta ,
\gamma \} \in [\Gamma ]^2$,
\sm
\ce{$T_\tau (\bar P^\lambda ) \forces_\MI " \mu (\breve O^{
\alpha (\Gamma^\beta \cup \Gamma^\gamma ) } \setminus
(H_\tau^{\alpha (\Gamma^\beta )} \cup H_\tau^{\alpha
(\Gamma^\gamma)} )) \leq {2 \over n^4} "$,}
\sm
\no hence by Lemma 5.8., $\mu (H_\tau^{\alpha (\Gamma^\beta
\cup \Gamma^\gamma )} \setminus (H_\tau^{\alpha (\Gamma^\beta)}
\cup H_\tau^{\alpha (\Gamma^\gamma)} )) \leq {2 \over n^4}$.
Therefore we get $\mu (B \setminus (H_\tau^{\alpha (\Gamma^\beta
)}
\cup H_\tau^{\alpha (\Gamma^\gamma)} )) \leq {3 \over n^4}$ for
all $\{ \beta , \gamma \} \in [\Gamma ]^2$. \par
Note, however, that
$$\eqalign{2 &\leq \sum_{\beta \in \Gamma} \mu (B \setminus
H_\tau^{\alpha (\Gamma^\beta)}) \cr &\leq \mu (\{ x;\; \vert
\{ \beta \in\Gamma ;\; x \in B \setminus H_\tau^{\alpha (\Gamma^\beta
)} \} \vert \leq 1 \}) + 2 \cdot n \cdot \mu ( \{ x ;\;
\vert \{ \beta\in\Gamma ;\; x \in B \setminus H_\tau^{\alpha
(\Gamma^\beta)} \} \vert \geq 2 \} ) \cr
&\leq 1 + 2 \cdot n \cdot n^2 \cdot {3 \over n^4} = 1
+ {6 \over n} < 2 ,\cr}$$
a contradiction to (***). Applying the partition relation
$\omega_1 \to (\omega , \omega_1 )^2$, we easily get
$\Delta_i \in [\Delta_{i-1}]^{\omega_1}$ so that
for all pairs $\{ \beta, \gamma \} \in [\Delta_i ]^2$ we have
$\mu (B \setminus H_{\sigma_i}^{\alpha (\Gamma^\beta \cup
\Gamma^\gamma)} ) > 0$. This proves the Main Lemma.
$\qed$
\bigskip
{\bolds 5.11. Observation.} {\it Let $G$ be $\QQ$--generic
over $N \prec H (\chi)$, where $\lambda \sub N$, $\QQ \in
N$, and $\vert N \vert = \lambda$.
\par
(i) $\hat S := \cup \{ \Sigma^p ;\; p\in G \}$ is the
set of splitting nodes of a Miller tree $S$; in fact,
whenever $\sigma \in \hat S$, then $\vert P^\lambda_\sigma
\cap \hat S \vert = \omega$.
\par
(ii) For $\sigma \in \hat S$ and $i\in\omega$ set $\hat
B_{\sigma , i} := \cap \{ B_{\sigma, i}^p ;\; p\in G\}$.
Then $\hat B_{\sigma , i}$ is a closed (and non--empty)
set of reals, $i < j$ implies $\hat B_{\sigma ,i} \supseteq
\hat B_{\sigma , j}$, $\tau \in P_\sigma \cap \hat S$ and
$\tau (lh (\sigma)) \geq i$ imply $\hat B_{\tau , j} \sub
\hat B_{\sigma , i}$ for any $j$, and for any $f \in [S]$,
we have $\vert \bigcap_{f \restrict n \in \hat S} \hat B_{f
\restrict n , 0} \vert = 1$ (and this is still true for
any branch of $S$ in any larger model; thus the family
$\{ \hat B_{\sigma , 0} ;\; \sigma\in \hat S\}$ can be
thought of as a $\MI$--name for a real).
\par
(iii) $\forall \beta < \lambda \; \exists n\in\omega\;
\forall\sigma\in\hat S\;\forall^\infty i\in\omega
\; (\hat B_{\sigma , i} \cap H_\sigma^{(\{ (\beta, n) \} )}
= 0)$.
}
\sm
{\capit Proof.} (i) Use genericity and the fact that given
$p, \sigma \in \Sigma^p$ and $i\in\omega$, we can find
find $q\leq p , \tau \in P_\sigma$ with $\tau (lh (\sigma))
\geq i$ and $\tau \in \Sigma^q$.
\par
(ii) Using genericity and the fact that given $p , B_{\sigma , i}^p$
we
can find $q \leq p$ so that $B^q_{\sigma , i}$ is closed,
we get closedness. The rest is easy.
\par
(iii) By genericity it suffices to show that given $p ,
\beta < \lambda$, there are $n \in\omega$ and $q \leq p$
with $(\beta , n ) \in \Gamma^q$.
To see this we use an argument similar to the one in the proof
of 5.10. Let $\Sigma^p = \{ \sigma_i ; \; i \in k \}$.
Using $\breve O^{\alpha (\Gamma^p)}$--maximality of $\bar P^\lambda
$, find $n_i$ ($i\in k$) so that $\forall^\infty \tau \in
P^\lambda_{\sigma_i} \; (\mu (B_i \setminus H_\tau^{\alpha (\Gamma^p
)}) \geq {1 \over n_i})$, where $B_i = B_{\sigma_i , n_{\sigma_i}
- 1}$ (this is (**) in 5.10.). Let $n > \max_{i\in k}
n_i$. We set $\Gamma^q = \Gamma^p \cup \{ (\beta , n )\}$.
Fix $i\in k$. As in 5.10. $\mu (H^{\alpha (\Gamma^q)}_\tau
\setminus (H_\tau^{\alpha (\Gamma^p)} \cup H_\tau^{\alpha
( \{ (\beta , n) \} )} )) \to 0$ for $\tau (lh (\sigma_i ))
\to \infty$ ($\tau\in P^\lambda_{\sigma_i}$). Also
(by Lemma 5.8.) $\mu (H^{\alpha ( \{ ( \beta , n ) \} ) }_\tau
) \leq {1 \over n}$. Hence $\liminf_{\tau\in P^\lambda_{\sigma_i}}
\mu (B_i \setminus H_\tau^{\alpha (\Gamma^q)}) \geq 
{1 \over n_i} - {1 \over n}$. Thus $\mu (B_i \setminus
H_{\sigma_i}^{\alpha (\Gamma^q)} ) > 0$. This means
we can make a straightforward extension to a condition
$q$. $\qed$
\bigskip
{\bolds 5.12. Claim.} {\it $S \forces_\MI " \forall
\alpha < \lambda \; (\breve x \not\in \breve N^\alpha )"$,
where $\breve x$ is the $\MI$--name given by 5.11. (ii).}
\sm
{\capit Proof.} Suppose, by way of contradiction, that
for some $\beta < \lambda$ and some $T \leq S$,
\sm
\ce{$T \forces_\MI " \breve x \in \breve N^\beta "$.}
\sm
\no Choose $n\in\omega$ so that for some $p \in G$,
$(\beta , n) \in \Gamma^p$. Clearly $T \forces_\MI
" \breve x \in \breve U^\beta_n = \breve O^\alpha "$,
where $\alpha = \alpha ( \{ (\beta,n) \} )$. Find
$k \in \omega$ and $T' \leq T$ so that
\sm
\ce{$T ' \forces_\MI " \breve x \in s_{\breve f^\alpha (k)}
"$.}
\sm
\no Let $\tau\in split (T')$ so that $lh (f^\alpha (\tau))
\geq k+1$. Then $T_\tau ' \forces_\MI " \breve x \in
s_{f^\alpha (\tau) (k)} \sub H^\alpha_\tau "$ (+).
Choose $q \leq p$, $q \in G$ so that $\tau\in\Sigma^q$.
Then $H^\alpha_\tau \cap B^q_{\tau , n_\tau^q - 1} = \emptyset$,
in particular $s_{f^\alpha (\tau) (k)} \cap
B^q_{\tau , n_\tau^q - 1} = \emptyset$. Choose $\tilde\tau
\in P^\lambda_\tau \cap T'$ with $\tilde\tau (lh (\tau))
\geq n^q_\tau - 1$, and $r \leq q$, $r \in G$
so that $\tilde\tau \in \Sigma^r$. 
Then $\hat B_{\tilde\tau , 0} \sub B^r_{\tilde\tau , 0}
\sub B^q_{\tau , n_\tau^r - 1}$, and
\sm
\ce{$T_{\tilde\tau} ' \forces_\MI " \breve x \in \hat
B_{\tilde\tau , 0} "$,}
\sm
\no contradicting (+). $\qed$
\sm
Thus $S$ and $\breve x$ satisfy $(*)$, and Theorem 5.6.
is proved. $\qed$
\bigskip
Theorems 5.1. and 5.6. say that the behaviour of
Miller forcing is rather similar to that of
Cohen forcing: both preserve $MA(\sigma$--centered)
and both collapse $add({\cal L})$ to $\omega_1$
(cf the first paragraph of the Introduction for
Cohen forcing). It seems that any "reasonable" 
forcing notion generated by a name for an unbounded
real does the latter (compare this to Question 3.11.).
\par
We believe that a similar argument as the one in part (d)
above should yield the same result for adding a Laver
real (and thus positively answer Question 3.7.).
The main problem is that when defining $\breve O$--softness
$(\heartsuit)$ we cannot require something corresponding
to condition (IV) for Laver forcing; and thus we do not know whether
the p.o. $\QQ$ will be $ccc$ in the Laver case.
\par
We also note that it seems to be a general state of affairs
that "reasonable" p.o.'s adding a single real can only
collapse cardinal invariants (in $ZFC$) which are dual
to cardinals which are increased by iterating the same 
forcing. E.g.: a single Cohen real collapses $cov ({\cal
L})$ and iterated Cohen forcing increases $unif ({\cal L})$;
or: a single Laver real collapses ${\bf d}$ and iterated
Laver forcing increases ${\bf b}$. According to this philosophy,
it should be consistent that $cov ({\cal L})$ is still large
after one Laver real, for Woodin (unpublished) and (later,
but independently) Judah and Shelah [JS 2, section 1]
showed that Laver forcing preserves outer measure
(and thus the iteration does not increase $unif ({\cal L})$.
\par
Of course one may investigate other related forcings
as well --- our choice was motivated by selecting forcings which
have some effect on cardinal invariants when iterated.
An example which leaves all cardinal invariants small
(and just increases $2^\omega$) is iterated Sacks forcing.
This seemingly corresponds to a result of Carlson's and
Laver's [CL] that Sacks forcing may preserve the full
extent of $MA$.

\Bigskip
\centerline{\capitg References}
\Smallskip
\itemitem{[Ba]} {\capit T. Bartoszy\'nski,} {\it Combinatorial
aspects of measure and category,} Fundamenta Mathematicae, vol. 127
(1987), pp. 225-239.
\smallskip
\itemitem{[BaJS]} {\capit T. Bartoszy\'nski, H. Judah and S.
Shelah,} {\it The Cicho\'n diagram,} submitted to Journal of Symbolic
Logic.
\smallskip
\itemitem{[BS]} {\capit T. Bartoszy\'nski and S. Shelah,}
{\it Closed measure zero sets,} Annals of Pure
and Applied Logic, vol. 58 (1992), pp. 93-110.
\smallskip
\itemitem{[Be]} {\capit M. G. Bell,} {\it On the combinatorial
principle $P(c)$,} Fundamenta Mathematicae, vol. 114 (1981),
pp. 149-157.
\smallskip
\itemitem{[Bl]} {\capit A. Blass,} {\it Selective ultrafilters 
and homogeneity,} Annals of Pure and Applied Logic,
vol. 38 (1988), pp. 215-255.
\smallskip
\itemitem{[Br 1]} {\capit J. Brendle,} {\it Larger cardinals
in Cicho\'n's diagram,} Journal of Symbolic Logic,
vol. 56 (1991), pp. 795-810.
\smallskip
\itemitem{[Br 2]} {\capit J. Brendle,} {\it Amoeba--absoluteness and
projective measurability,} to appear in Journal of Symbolic Logic.
\smallskip
\itemitem{[BJS]} {\capit J. Brendle, H. Judah and S. Shelah,}
{\it Combinatorial properties of Hechler forcing,} Annals
of Pure and Applied Logic, vol. 58 (1992), pp. 185-199.
\smallskip
\itemitem{[CL]} {\capit T. Carlson and R. Laver,} {\it Sacks
reals and Martin's axiom,} Fundamenta Mathematicae,
vol. 133 (1989), pp. 161-169.
\smallskip
\itemitem{[CP]} {\capit J. Cicho\'n and J. Pawlikowski,} {\it
On ideals of subsets of the plane and on Cohen reals,}
Journal of Symbolic Logic, vol. 51 (1986), pp. 560-569.
\smallskip
\itemitem{[Fr]} {\capit D. Fremlin,} {\it Cicho\'n's diagram,}
S\'eminaire Initiation \`a l'Analyse (G. Choquet,
M. Rogalski, J. Saint Raymond), Publications Math\'ematiques
de l'Universit\'e Pierre et Marie Curie, Paris, 1984,
pp. 5-01 - 5-13. 
\smallskip
\itemitem{[GJS]} {\capit M. Goldstern, M. Johnson and O. Spinas,}
{\it Towers on trees,} to appear in Proceedings of the American
Mathematical Society.
\smallskip
\itemitem{[GRSS]} {\capit M. Goldstern, M. Repick\'y, S.
Shelah and O. Spinas,} {\it On tree ideals,} preprint.
\smallskip
\itemitem{[Ha]} {\capit L. Halbeisen,} in preparation.
\smallskip
\itemitem{[Je 1]} {\capit T. Jech,} {\it Set theory,} Academic Press,
San Diego, 1978.
\smallskip
\itemitem{[Je 2]} {\capit T. Jech,} {\it Multiple forcing,}
Cambridge University Press, Cambridge, 1986.
\smallskip
\itemitem{[Ju 1]} {\capit H. Judah,} {\it $\Sigma^1_2$--sets
of reals,} Journal of Symbolic Logic, vol. 53 (1988), pp. 636-642.
\smallskip
\itemitem{[Ju 2]} {\capit H. Judah,} {\it Absoluteness for projective
sets,} to appear in Logic Colloquium 1990.
\smallskip
\itemitem{[JMS]} {\capit H. Judah, A. Miller and S. Shelah,}
{\it Sacks forcing, Laver forcing, and Martin's Axiom,}
Archive for Mathematical Logic, vol. 31 (1992), pp. 145-162.
\smallskip
\itemitem{[JS 1]} {\capit H. Judah and S. Shelah,} {\it 
$\Delta_2^1$--sets of reals,} Annals of Pure and Applied Logic,
vol. 42 (1989), pp. 207-223.
\smallskip
\itemitem{[JS 2]} {\capit H. Judah and S. Shelah,} {\it The Kunen-Miller
chart (Lebesgue measure, the Baire property, Laver reals and
preservation theorems for forcing),} Journal of
Symbolic Logic, vol. 55 (1990), pp. 909-927.
\smallskip
\itemitem{[Ku]} {\capit K. Kunen,} {\it Set theory,} North-Holland,
Amsterdam, 1980.
\smallskip
\itemitem{[Ma]} {\capit A. R. D. Mathias,} {\it Happy families,}
Annals of Mathematical Logic, vol. 12 (1977), pp. 59-111.
\smallskip
\itemitem{[Mi 1]} {\capit A. Miller,} {\it Some properties of
measure and category,} Transactions of the American Mathematical
Society, vol. 266 (1981), pp. 93-114.
\smallskip
\itemitem{[Mi 2]} {\capit A. Miller,} {\it Rational perfect
set forcing,} Contemporary Mathematics, vol. 31 (Axiomatic
Set Theory, 1984, edited by J. Baumgartner, D. Martin and
S. Shelah), pp. 143-159.
\smallskip
\itemitem{[Pa 1]} {\capit J. Pawlikowski,} {\it Powers of
transitive bases of measure and category,} Proceedings
of the American Mathematical Society, vol. 93 (1985),
pp. 719-729.
\smallskip
\itemitem{[Pa 2]} {\capit J. Pawlikowski,} {\it Why Solovay
real produces Cohen real,} Journal of Symbolic Logic,
vol. 51 (1986), pp. 957-968.
\smallskip
\itemitem{[Ro]} {\capit J. Roitman,} {\it Adding a random or a
Cohen real: topological consequences and the effect on
Martin's axiom,} Fundamenta Mathematicae, vol. 103 (1979),
pp. 47-60 and vol. 129 (1988), p. 141.
\smallskip
\itemitem{[RS]} {\capit A. Ros{\l}anowski and S. Shelah,}
{\it More forcing notions imply diamond,} preprint.
\smallskip
\itemitem{[SW]} {\capit S. Shelah and H. Woodin,}
{\it Forcing the failure of CH by adding a real,}
Journal of Symbolic Logic, vol. 49 (1984), pp. 1185-1189.
\smallskip
\itemitem{[Si]} {\capit J. Silver,} {\it Every analytic set is
Ramsey,} Journal of Symbolic Logic, vol. 35 (1970),
pp. 60-64.
\smallskip
\itemitem{[Tr 1]} {\capit J. Truss,} {\it Sets having calibre
$\aleph_1$,} in: R. O. Gandy and J. M. E. Hyland, eds.,
Logic Colloquium '76 (North--Holland, Amsterdam, 1977),
pp. 595-612.
\smallskip
\itemitem{[Tr 2]} {\capit J. Truss,} {\it The noncommutativity
of random and generic extensions,} Journal of Symbolic Logic,
vol. 48 (1983), pp. 1008-1012.
\smallskip
\itemitem{[Tr 3]} {\capit J. Truss,} {\it Connections between
different amoeba algebras,} Fundamenta Mathematicae, vol.
130 (1988), pp. 137-155.
\smallskip
\itemitem{[vD]} {\capit E. K. van Douwen,} {\it
The integers and topology,} Handbook of set--theoretic
topology, K. Kunen and J. E. Vaughan (editors),
North--Holland, Amsterdam, 1984, pp. 111-167.
\smallskip
\itemitem{[Wo]} {\capit W. H. Woodin,} {\it On the consistency strength
of projective uniformization,} Proceedings of the Herbrand Symposium,
Logic Colloquium '81, J. Stern (ed.), North--Holland, Amsterdam,
1982, 365-384.
\smallskip

\vfill\eject\end